\documentclass[moor,nonblindrev]{informs1} % for review, not blinded
\usepackage{graphicx, subfig, float, mathtools}
\usepackage{algorithm}
\usepackage{algorithmic}
\usepackage{listings}
\usepackage{enumerate}
\usepackage{color}%Liu Hongying 2020.08.03

\usepackage{natbib}
 \NatBibNumeric
 \def\bibfont{\small}%
 \bibpunct[, ]{[}{]}{,}{n}{}{,}%

\usepackage[colorlinks=true,breaklinks=true,bookmarks=true,urlcolor=blue,
     citecolor=blue,linkcolor=blue,bookmarksopen=false,draft=false]{hyperref}

\def\EMAIL#1{\href{mailto:#1}{#1}}% When hyperref is used, otherwise outcomment
\def\URL#1{\href{#1}{#1}}         % When hyperref is used, otherwise outcomment

\def\eproof{\hspace*{\fill}~{\setlength{\fboxsep}{0pt}\setlength{\fboxrule}{0.2pt}\fbox{\rule[0pt]{0pt}{1.3ex}\rule[0pt]{1.3ex}{0pt}}}\\}%

\def\ST{\songti\rm\relax}

\def\ST{{\rm s.t.}}

\def\R{{\mathbb R}}%\def\R{\mbox{I}\!\mbox{R}}

\newtheorem{que}{Question}[section]
\TheoremsNumberedThrough     % Preferred (Theorem 1, Lemma 1, Theorem 2)
\EquationsNumberedThrough    % Default: (1), (2), ...
\begin{document}
%%%%%%%%%%%%%%%%

%\setcounter{page}{1}
%\DOI{10.1287/xxxx.0000.0000}%

\RUNAUTHOR{M. Song et al.}

 \RUNTITLE{Local Optimality for Hidden Convex Optimization}
%\RUNTITLE{}

\TITLE{Local Optimality Conditions for a Class of Hidden Convex Optimization}
%\TITLE{Global and Local Optimality Conditions for a Class of Hidden Optimization Combining Trust-region Subproblem with Convex Optimization}

\ARTICLEAUTHORS{
\AUTHOR{Mengmeng Song$^1$}
%\AFF{$1$, School of Mathematical Sciences, Beihang University, Beijing, 100191, People's Republic of China,
\EMAIL{songmengmeng@buaa.edu.cn}, \URL{}
%}
\AUTHOR{Yong Xia$^1$}
%\AFF{School of Mathematical Sciences, Beihang University, Beijing, 100191, People's Republic of China,
\EMAIL{yxia@buaa.edu.cn}, \URL{}
%}
\AUTHOR{Hongying Liu$^1$, Corresponding Author}
\AFF{\EMAIL{liuhongying@buaa.edu.cn}, \URL{}\\
	$1$, School of Mathematical Sciences, Beihang University, Beijing, 100191, People's Republic of China,
}
} % end of the block
\ABSTRACT{
Hidden convex optimization is such a class of nonconvex optimization problems that can be globally solved in polynomial time via equivalent convex programming reformulations.  In this paper, we focus on checking local optimality in hidden convex optimization.  We first introduce a class of hidden convex optimization problems by jointing the classical nonconvex trust-region subproblem (TRS) with convex optimization (CO), and then present a comprehensive study on local optimality conditions. In order to guarantee the existence of a necessary and sufficient condition for local optimality, we need more restrictive assumptions. To our surprise, while (TRS) has at most one local non-global minimizer and (CO) has no local non-global minimizer, their joint problem could have more than one local non-global minimizer.
}

%\KEYWORDS{deterministic inventory theory; infinite linear programming duality;
%  existence of optimal policies; semi-Markov decision process; cyclic schedule}
%\MSCCLASS{Primary: 90B05; secondary: 90C40, 90C90}
%\ORMSCLASS{Primary: Inventory/production: deterministic multi-item;
%  secondary: dynamic programming/optimal control: deterministic
%  semi-Markov; programming: infinite dimensional}
%\HISTORY{Received November 20, 2003; revised March 8, 2004, and March 26, 2004.}

\MSCCLASS{ 90C26,  90C46, 90C30}
%\ORMSCLASS{Primary: ; secondary: }
\HISTORY{}

\maketitle
%%%%%%%%%%%%%%%%%%%%%%%%%%%%%%%%%%%%%%%%%%%%%%%%%%%%%%%%%%%%%%%%%%%%%%
% Samples of sectioning (and labeling) in MOOR.
% NOTE: (1) all section levels end with a period,
%       (2) capitalization is as shown (sentence style, not title style).
%
%\section{Introduction.}\label{intro} %%1.
%\subsection{Duality and the classical EOQ problem.}\label{class-EOQ} %% 1.1.
%\subsection{Outline.}\label{outline1} %% 1.2.
%\subsubsection{Cyclic schedules for the general deterministic SMDP.}
%  \label{cyclic-schedules} %% 1.2.1
%\section{Problem description.}\label{problemdescription} %% 2.
% Text of your paper here

\section{Introduction}
%介绍部分首先从凸优化开始（凸优化广受关注，引一些重要参考书，特点凸优化局部解就
%是全局），一般非凸优化太困难，研究一个“距离凸优化不太远”的非凸优化，注意到TRS 是一个well-studied 非凸优化，因此引入这一个联合非凸优化问题。 这个问题还包含一些特殊应用，比如p- 子问题。
As a fundamental topic in operations research,
convex optimization has a rapid development in the last two decades \cite{boyd2004convex,Nesterov2018}, thanks to   its popular applications in machine learning \cite{bubeck2015convex}.  %mathematical optimization problems, includes least-squares, linear programming problems, semidefinite programs and second-order cone programs.
%The development of convex optimization in the last two decades has been very rapid and exciting, see, for example,   \cite{boyd2004convex,Nesterov2018}.
One of the key to the algorithmic success in convex optimization is that any local minimizer is in fact a global one. The situation is different in nonconvex optimization, which is in general NP-hard, and often has many local non-global minimizers. Local optimal solutions play a great role in globally solving structured nonconvex optimization problem \cite{B17}. However, as shown in \cite{M87,P88}, checking local optimality of a nonconvex quadratic program is already NP-hard. That is, local optimization is as difficult as global optimization. On the other hand, we notice that not all nonconvex optimization problems are difficult to globally solve. For example, hidden convex optimization admits an equivalent polynomial-solvable convex programming reformulation, see \cite{X20} for a recent survey.  
Accordingly, we ask
\begin{que}\label{ques}
~Can we go far beyond the standard local optimality conditions for hidden convex optimization?
\end{que}

A typical hidden convex optimization is the following trust region subproblem (TRS):
\begin{equation}\label{eq:TRS}\tag{TRS}
\begin{array}{cl}
\min\limits_{x\in \R^n} \{\tfrac{1}{2}x^THx+c^Tx: x^Tx\le \Delta\},
\end{array}
\end{equation}
where $H=H^T\in \R^{n\times n}, c\in \R^n$ and $\Delta>0$. (TRS) plays a key role in trust region methods for solving nonlinear programming problems \cite{Conn00, yuan15}.
%It's somewhat of a surprise that there is no gap between \eqref{eq:TRS} and its Lagrangian dual problem, i.e., strong duality holds. Actually,
In the early 1980s, Gay \cite{Gay81}, Sorensen \cite{Sorensen82}, Mor\'{e} and Sorensen \cite{More83} have established necessary and sufficient optimality condition for the global minimizer of \eqref{eq:TRS}. Observing this, Ye \cite{Y92} showed that (TRS) can be globally solved in polynomial time.  In 1994, Mart\'{i}nez  \cite{JOSE94} proved that \eqref{eq:TRS} has at most one local non-global minimizer, and then established a necessary optimality condition and a  sufficient one for the local non-global minimizer of \eqref{eq:TRS}. Observing the gap between necessary and sufficient conditions, Tao and An \cite{Tao98}
pointed out that ``For the trust-region subproblem, it is interesting to note the following paradox: checking a global solution is easier than checking a local nonglobal one.''  It seems that the answer to Question \ref{ques} is negative on (TRS). %Actually, it followed by the conclusion in \cite{Jiulin20} that there is necessary and sufficient consition for local nonglobal minimizer of \eqref{eq:TRS}. So, by now we may assert that it is equally easy for checking the global minimizer and local nonglobal minimizer. That breaks the paradox proposed in \cite{Tao98}.
In 2020, Wang and Xia \cite{Jiulin20} broke Tao and An's paradox by proving that Mart\'{i}nez's sufficient optimality condition \cite{JOSE94} for \eqref{eq:TRS} is also necessary. In the same paper, it is also shown that the local non-global minimizer can be found or proved to do not exist in polynomial time. Thus, Wang and Xia's result gives a positive answer to Question \ref{ques} on \eqref{eq:TRS}.

%Since a general nonconvex optimization is too hard to deal with, we will study a one which is near by convex problem.

Towards answering Question \ref{ques},  we introduce in this paper a class of nonconvex optimization, which joints nonconvex trust-region subproblem with  convex optimization:
\begin{equation}\label{eq:TRS-C}\tag{TRS-C}
\begin{array}{cl}
\min\limits_{x\in \R^n, y\in\R^m} & \tfrac{1}{2}x^THx+c^Tx+f_0(y)\\
\ST   &x^Tx+ f_1(y)\le 0\\
 &\ \ \ \ \ \ \ \ \, f_j(y)\le 0, j=2,\cdots,k,
\end{array}
\end{equation}
where $H=H^T\in \R^{n\times n}, c\in \R^n, k\in \mathbb{Z}^+$, and
$f_j(\cdot)~ (j=0, 1,\cdots,k)$ are convex and twice continuously differentiable. Throughout this paper, we assume $H\nsucceq 0$, otherwise \eqref{eq:TRS-C} is a convex optimization problem. In particular, setting $m=1, k=1, f_0\equiv0, f_1\equiv-\Delta$, \eqref{eq:TRS-C} reduces to \eqref{eq:TRS}.

As we have mentioned above, \eqref{eq:TRS} has at most one local non-global minimizer, and convex optimization problem has no local non-global minimizer. %It is natural to extend this property to  \eqref{TRS-C} as it joints \eqref{eq:TRS} with  convex optimization. Namely, we ask
Notice that \eqref{eq:TRS-C} joints \eqref{eq:TRS} with  convex optimization. So, it is natural to expect a positive answer to the following question:
\begin{que}\label{ques2}
~Does \eqref{eq:TRS-C} have at most one local non-global minimizer?
\end{que}

%When $H=0, c=0$, \eqref{eq:TRS-C} is
%\begin{equation*}
%\min_{y\in \R^m} \{f_1(y): f_j(y)\le 0, j=2,\cdots,k\}.
%\end{equation*}
%which is a general convex optimization problem.
%Particularly, both of the well-known TRS and the $p-$regularized subproblem are the instance of \eqref{eq:TRS-C}.

We may consider another special case of \eqref{eq:TRS-C} by
setting $m=1, k=1, f_0(y)=\tfrac{\sigma}{p}y^{\tfrac{p}{2}}, f_1(y)=-y$, where $\sigma>0$ and $p>2$ are two parameters:
\begin{equation}\label{eq:pRS-gtrs}
\min\limits_{x\in \R^n,y\in\R}\{ \tfrac{1}{2}x^THx+c^Tx+\tfrac{\sigma}{p}y^{\tfrac{p}{2}}: x^Tx-y\le 0\}.
\end{equation}
The univariate variable $y$ can be removed by substituting the constraint into the objection. Then \eqref{eq:pRS-gtrs} is equivalent to the following unconstrained optimization
  \begin{equation}\label{eq:pRS}\tag{$p$-RS}
\min\limits_{x\in \R^n}  \tfrac{1}{2}x^THx+c^Tx+\tfrac{\sigma}{p}\|x\|^p,
\end{equation}
which is known as $p-$regularized subproblem \cite{Gould10}.
%arising in the $p-$regularized methods \cite{Griewank81}.
%Thus, the general problem \eqref{eq:TRS-C} we consider is extension of convex optimization problem, trust region subproblem and $p$-regularized subproblem.
%In literature, the most common choice to regularize the quadratic approximation is
In particular, the cubic regularization ($p=3$) was first introduced in \cite{Griewank81} and then widely applied in nonlinear programming, see for example, \cite{Cartis11, Y.Nesterov06, Weiser07}. The other special case of \eqref{eq:pRS} with $p=4$ corresponds to the double-well potential optimization \cite{F17,XS17}.
%Recently, a comprehensive comparison for the numerical effectiveness \textcolor{red}{between \eqref{eq:pRS} for general $p>2$ and \eqref{eq:TRS} was made in \cite{Gould10}.
Recently, Hsia et al. \cite{Xia17} established necessary and sufficient optimality conditions for both global and local minimizers of  \eqref{eq:pRS} with any $p>2$.  In the same paper, \eqref{eq:pRS} is shown to have at most one local non-global minimizer. Thus, both Question \ref{ques} and \ref{ques2} have positive answers on \eqref{eq:pRS}.

%The local non-global minimizer of \eqref{eq:TRS-C} plays a great role in globally solving the extended \eqref{eq:TRS-C}, which has more additional constraints to \eqref{eq:TRS-C}:
%\begin{equation}\label{eq:ETRS-C}
%\begin{array}{cl}
%\min\limits_{x\in \R^n, y\in\R^m} & \tfrac{1}{2}x^THx+c^Tx+f_0(y)\\
%\ST   & x^Tx+f_1(y)\le 0\\
%& f_j(y)\le 0, j=2,\cdots,k\\
%& g(y)\le 0
%\end{array}
%\end{equation}
%where $g:\R^m\rightarrow\R$ is continuous. Let $X^*$ be the set of globally optimal solutions of \eqref{eq:TRS-C} and $X^+$ be the set of local non-global optimal solutions of \eqref{eq:TRS-C}. Consider the nontrivial case $X^*\cap\{(x,y):g(y)\le 0\}=\emptyset$. Then, the optimal solution of \eqref{eq:ETRS-C} is either
%in $X^+$ or the optimal solution of the following equality version:
%\begin{equation*}
%\begin{array}{cl}
%\min\limits_{x\in \R^n, y\in\R^m} & \tfrac{1}{2}x^THx+c^Tx+f_0(y)\\
%\ST   & x^Tx+f_1(y)\le 0\\
%& f_j(y)\le 0, j=2,\cdots,k\\
%& g(y)= 0.
%\end{array}
%\end{equation*}

In this paper, we give a comprehensive study on \eqref{eq:TRS-C}. We first reveal its hidden convexity by establishing necessary and sufficient optimality condition for global minimizer. It makes sense as \eqref{eq:TRS-C} joints a hidden convex with a convex optimization problem. We then establish necessary and sufficient optimality conditions for local non-global minimizer, respectively. However, the gap of two conditions exists for the general case. In order to close this gap, we need more assumptions. It supports the conjecture that Question \ref{ques} has no positive answer. For Question \ref{ques2}, to our surprise,  the answer is negative. Actually,    \eqref{eq:TRS-C} may have a finite number of local non-global minimizers.  Only for very special cases including \eqref{eq:TRS} and \eqref{eq:pRS},
\eqref{eq:TRS-C} has at most one local non-global minimizer.
Consequently,  we conclude that {\it
the existence of many local minimizers is NOT a (or at least not a unique) reason making
global optimization difficult.
}

%Introduction 再加一个surprising 对于(TRS-C）， 刻画局部解要比刻画全局解远远困难。

The remainder of this paper is organized as follows.
Section \ref{sec:3} presents necessary and sufficient conditions for global minimizers of \eqref{eq:TRS-C}.  Section \ref{sec:4} establishes  necessary and  sufficient conditions  for local non-global minimizer of \eqref{eq:TRS-C}, respectively. These conditions are further characterized  by a scalar function
for a special case of \eqref{eq:TRS-C} with $k=1$ in Section \ref{sec:5}.
We consider another special case of \eqref{eq:TRS-C} where $f_1(y)$ is a linear function in Section \ref{sec:6}.  Two instances with more than one local non-global minimizer are constructed. Then a sufficient condition is presented to guarantee that \eqref{eq:TRS-C} has at most one local non-global minimizer. As an extension of \eqref{eq:TRS} and \eqref{eq:pRS}, we identify a class of cases of \eqref{eq:TRS-C} where local non-global minimizer enjoys a necessary and sufficient optimality condition. We conclude the paper
in the last section with a few open questions.

\emph{\bf Notations.}
For any matrix $P\in\R^{n\times n}$, $P^T$ and $\det{(P)}$ denote the transposition and determination of $P$, respectively.
$P\succ(\succeq)0$ denotes that $P$ is positive (semi)definite.
Let $I$ stand for the identity matrix of order $n$.
For a vector $x\in\R^n$, ${\rm Diag}(x)$ returns a diagonal matrix with diagonal elements being $x_1,\cdots,x_n$.
%For the sake of simplicity, we denote $q(x):=\tfrac{1}{2}x^THx+c^Tx$ throughout this paper.
 % $e$ denotes a vector of dimension $n$ with all components equal to one. For a number $\beta\in\R$, ${\rm sign}(\beta)=\beta/|\beta|$, if $\beta\neq 0$, otherwise ${\rm sign}(\beta)=0$.
Let $f^{-1}(\cdot)$ represent the inverse of function $f(\cdot)$, if it exists. $\mathbb{Z}^+$ denotes the set of all positive integers.
The eigenvalue decomposition of the symmetric matrix $H$ is given by
\begin{equation}\label{eq:tzzfjH}
H=V{\rm Diag}(\lambda_1,\cdots,\lambda_n)V^T,
\end{equation}
where $V=(v_1,v_2,\cdots,v_n)\in\R^{n\times n}$ is orthogonal, and %$\lambda_1\le\lambda_2\le\cdots \le\lambda_n$.
$\lambda_i$ is the $i$-th smallest eigenvalue of $H$.

\section{Global Optimality Condition}\label{sec:3}
In this section, we present the necessary and sufficient optimality condition for the global minimizer of \eqref{eq:TRS-C}, where the Slater condition is assumed.

\begin{theorem}\label{theo:qjjbytj}
Assume that $(x_*,y_*)\in \R^n\times \R^m$ is a global minimizer of \eqref{eq:TRS-C} and there exists $\bar y$ such that $f_j(\bar y)<0,~j=1,\cdots,k$. There exist $\mu^*_j\ge 0$ ($j=1,\cdots,k$) such that
\begin{eqnarray}
&&(H+\mu^*_1 I)x_*+c=0,\label{eq:Lx}\\
&&\nabla f_0(y_*)+\frac{\mu^*_1}{2}\nabla f_1(y_*)+\sum_{j=2}^{k}\mu^*_j\nabla f_j(y_*)=0,\label{eq:Ly}\\
&&
x_*^Tx_*+f_1(y_*)=0,\ \  \mu^*_j f_j(y_*)=0, ~ j=2,\cdots,k,\label{eq:Lc}\\
&& H+\mu^*_1 I\succeq 0.\label{eq:Lmu}
\end{eqnarray}
%and $H+\mu^*_1 I$ is positive semidefinite.
\end{theorem}

\emph{Proof.}
%Let $\lambda_1$ be the smallest eigenvalue of $H$. If $\lambda_1\ge 0$, then $H$ is positive semidefinite and \eqref{eq:TRS-C} is a convex optimization problem. Since $(0,\bar y)$ is strictly feasible for \eqref{eq:TRS-C}, the Slater condition holds. According to first order necessary optimality conditions, there exist Lagrange multipliers $\mu^*_j\ge 0, j=1,\cdots,k$ satisfies \eqref{eq:Lx}-\eqref{eq:Lc}.
%Since $\mu^*_1\ge 0$ and $H$ is positive semidefinite, $H+\mu^*_1I$ is positive semidefinite.
%Now we assume $\lambda_1<0$ and hence $H$ is not positive semidefinite.
Note that $\lambda_1<0$ for $H\nsucceq 0$.
Since $(x_*,y_*)$ is a global minimizer of \eqref{eq:TRS-C}, then $x_*$ solves
\begin{equation}\label{eq:GTRSx}
\min_{x\in \R^n} \{\tfrac{1}{2}x^THx+c^Tx:~ x^Tx\le -f_1(y_*)\}.
\end{equation}
We conclude that
\begin{equation}\label{eq:qjjds}
x_*^Tx_*+f_1(y_*)=0.
\end{equation}
Otherwise, it follows from $x_*^Tx_*+f_1(y_*)<0$ that $H$ is positive semidefinite by second order necessary optimality, which contradicts the fact $\lambda_1<0$.
% for \eqref{eq:GTRSx}

Therefore, $(x_*,y_*)$ remains a global minimizer of
\begin{equation*}
\begin{array}{cl}
\min\limits_{x\in \R^n, y\in\R^m} & F(x,y):= \tfrac{1}{2}x^THx+c^Tx+f_0(y)-\frac{\lambda_1}{2}(x^Tx+f_1(y))\\
\ST   & x^Tx+f_1(y)=0\\
      & \ \ \ \ \ \ \ \,\  f_j(y)\le 0, ~j=2,\cdots,k.
\end{array}
\end{equation*}
%Furthermore, it can be shown that $(x_*,y_*)$ is a global minimizer of the convex program
Now, consider a convex optimization relaxation of the above problem:
\begin{equation}\label{eq:qjjdjtwt}
\begin{array}{cl}
\min\limits_{x\in \R^n, y\in\R^m} & F(x,y)=\tfrac{1}{2}x^T(H-\lambda_1 I)x+c^Tx+f_0(y)-\frac{\lambda_1}{2}f_1(y)\\
\ST   & x^Tx+f_1(y)\le 0\\
      & \ \ \ \ \ \ \ \ \, f_j(y)\le 0, j=2,\cdots,k.
\end{array}
\end{equation}
Let $(x,y)$ be any feasible solution to \eqref{eq:qjjdjtwt} with $x^Tx+f_1(y)<0$. Let $v_1$ be an eigenvector of $H$ corresponding to $\lambda_1$ (i.e., $Hv_1=\lambda_1v_1$) and satisfy $c^Tv_1\le 0$. Then there is a $\tau>0$ such that
$\bar x=x+\tau v_1$ satisfies $\bar x^T\bar x+f_2(y)=0$. We can verify that $(\bar x,y)$ is feasible to \eqref{eq:qjjdjtwt}, and
\begin{eqnarray}
\frac{1}{2}{\bar x}^T(H-\lambda_1I)\bar x+c^T\bar x
&=&\frac{1}{2}x^T(H-\lambda_1I) x+c^T x+\tau c^Tv_1\nonumber\\
&\le& \frac{1}{2}x^T(H-\lambda_1I) x+c^Tx\nonumber.
\end{eqnarray}
Therefore, $F(\bar x, y)\le F(x,y)$. It follows that $(x_*,y_*)$ must be a global minimizer of \eqref{eq:qjjdjtwt}.  Notice that \eqref{eq:qjjdjtwt} is a convex optimization problem where Slater condition holds by assumption. According to the first order necessary optimality condition, there exist $\bar\mu_1\ge 0$ and $\mu^*_j\ge 0, j=2,\cdots,k$, which satisfy
\begin{eqnarray}
&&(H+\bar\mu_1 I-\lambda_1 I)x_*=-c,\ \ \nabla f_0(y_*) + \frac{\bar\mu_1-\lambda_1}{2}\nabla f_1(y_*)+ \sum_{j=2}^{k}\mu^*_j\nabla f_j(y_*) =0,\label{eq:qjjdjKKT1}\\
&&
\bar\mu_1(x_*^Tx_*+f_1(y_*))=0,\ \ \mu^*_j f_j(y_*)=0, j=2,\cdots,k.\label{eq:qjjdjKKT2}
\end{eqnarray}
Define $\mu^*_1=\bar\mu_1-\lambda_1$.
Equations  \eqref{eq:qjjds}, \eqref{eq:qjjdjKKT1} and \eqref{eq:qjjdjKKT2} imply that $\mu^*_j, j=1,\cdots,k$, satisfy
\eqref{eq:Lx}, \eqref{eq:Ly} and \eqref{eq:Lc}.
Since $H-\lambda_1 I\succeq 0$ and $\bar\mu_1\ge 0$, we have $H+\mu^*_1 I\succeq 0$, i.e., \eqref{eq:Lmu} holds true.
\eproof

\begin{remark}
\eqref{eq:TRS-C} admits an equivalent convex reformulation
\eqref{eq:qjjdjtwt}, which reveals the hidden convexity of \eqref{eq:TRS-C}.
\end{remark}

Theorem \ref{theo:qjjbytj} establishes a necessary  condition for global optimality of $(x_*,y_*)$. The next result shows  that it is also a sufficient condition.

\begin{theorem}\label{theo:qjjcftj}
Let $(x_*,y_*)\in \R^n\times \R^m$ be a feasible solution to \eqref{eq:TRS-C}. If there exist $\mu^*_j\ge 0,  j=1,\cdots,k$, satisfying \eqref{eq:Lx}-\eqref{eq:Lmu}. Then $(x_*,y_*)$ is a global minimizer of \eqref{eq:TRS-C}.
\end{theorem}

\emph{Proof.}
%With the given conditions, $(x_*,y_*)$ is a KKT point of \eqref{eq:TRS-C}. If $\lambda_1\ge 0$,  \eqref{eq:TRS-C} belongs to convex optimization. Then, $(x_*,y_*)$ is a global minimizer of \eqref{eq:TRS-C}.
%Now we consider
Define $\bar \mu_1=\mu^*_1+\lambda_1$. It follows from $H+\mu^*_1 I\succeq 0$ that $\bar\mu_1\ge 0$. Substituting $\mu^*_1=\bar\mu_1-\lambda_1$ into \eqref{eq:Lx}, \eqref{eq:Ly} and \eqref{eq:Lc} yields \eqref{eq:qjjdjKKT1} and \eqref{eq:qjjdjKKT2}, respectively. That is, $(x_*,y_*)$ is a global minimizer of \eqref{eq:qjjdjtwt}. Therefore, for any feasible solution $(x,y)$ to \eqref{eq:TRS-C}, we have
$$
\tfrac{1}{2}x^THx+c^Tx+f_0(y)\ge F(x,y)\ge F(x_*,y_*)=\tfrac{1}{2}x_*^THx_*+c^Tx_*+f_0(y_*),
$$
where the first inequality follows from the feasibility of $(x,y)$ and the fact $\lambda_1<0$, the second inequality holds since $(x_*,y_*)$ is a global minimizer of \eqref{eq:qjjdjtwt}, and the equality is implied from \eqref{eq:Lc}. Thus,   $(x_*,y_*)$ is a global minimizer of \eqref{eq:TRS-C}.
\eproof

\section{Local Non-Global Optimality Conditions: General Case}\label{sec:4}

In this section, we present local non-global optimality conditions for \eqref{eq:TRS-C}.
%, where  $\lambda_1<0$ is assumed, since otherwise \eqref{eq:TRS-C} has no local non-global minimizer.
In contrast with linear independence constraint qualification (LICQ) made in general nonlinear programming, we need the following relaxed constraint qualification throughout this section:
\begin{assumption}\label{as:LICQ}
Let $(x_*,y_*)$ be any candidate local minimizer of \eqref{eq:TRS-C}. The gradients $\nabla f_i(y_*)$ for $i\in\{ j\in\{ 1,\cdots,k\}:~ f_j(y_*)=0\}$ are linearly independent.
\end{assumption}

%If $(x_*,y_*)$ is a local non-global minimizer of \eqref{eq:TRS-C}, then $x_*\neq 0$ and $x_*$ solves \eqref{eq:GTRSx}. What's more, the main conclusion  in Theorem \ref{th:local_pneq0} show that $\lambda_1<\lambda_2$, \eqref{eq:Lx},\eqref{eq:Ly} and \eqref{eq:Lc} hold at $(x_*,y_*)$ with $\max\{0,-\lambda_2\}<\mu^*_1<-\lambda_1$, $\mu^*_j\ge 0, j=2,\cdots,k$.

%In characterizing global minimizer, Slater condition can guarantee Lagrange multipliers existing since the hidden convex nature of \eqref{eq:TRS-C}. When characterizing local minimizer, the following assumption is proposed since the nonconvexity of \eqref{eq:TRS-C}.
%\begin{assumption}\label{as:LICQ}
%Let $(x_*,y_*)$ is feasible for \eqref{eq:TRS-C}. The gradients $\nabla f_i(y_*)$ for $i=1,\cdots,k$ with $f_i(y_*)=0$ are linearly independent.
%\end{assumption}
%\begin{remark}	
%In characterizing global minimizer, Slater condition is used to guarantee the existence of Lagrange multipliers. One can verify that, Theorem \ref{theo:qjjbytj} holds true,
%if we replace the Slater assumption ``there exists $\bar y$ such that $f_j(\bar y)<0, j=1,\cdots,k$'' made in Theorem  \ref{theo:qjjbytj} with Assumption \ref{as:LICQ}.
%\end{remark}
	
%First, the following lemma comes from traditional the first order and the second order necessary conditions and sufficient conditions for \eqref{eq:TRS-C}.
Applying the classical optimality conditions in nonlinear programming  to \eqref{eq:TRS-C} under Assumption \ref{as:LICQ}, we immediately have
\begin{lemma}\label{le:traoptcon}
\begin{enumerate}[(a)]
		\item
Let $(x_*,y_*)$ be a local minimizer of \eqref{eq:TRS-C}.
Under Assumption \ref{as:LICQ}, there exist $\mu^*_j\ge 0$ ($j=1,\cdots,k$) satisfying
\eqref{eq:Lx}, \eqref{eq:Ly}, and \eqref{eq:Lc}. Moreover,
\begin{eqnarray}
&\bmatrix
s^T & t^T
\endbmatrix&
\bmatrix
H+\mu^*_1 I & 0 \\
0 & \nabla^2 f_0(y_*)+\frac{\mu^*_1}{2} \nabla^2 f_1(y_*)+\sum_{j=2}^k\mu^*_j \nabla^2f_j(y_*)\\
\endbmatrix
\bmatrix
s \\
t \\
\endbmatrix
\ge 0\label{eq:jbejbytj}
\end{eqnarray}
holds for any $(s, t)\in T,$ where
\begin{eqnarray}\label{st}
&T=\{(s,t):&~2s^Tx_*+t^T\nabla f_1(y_*)=0;\\ \nonumber
&&~t^T\nabla f_i(y_*)=0, ~i\in\{j\in\{2,\cdots,k\}: f_j(y_*)=0, \mu^*_j>0\};\\ \nonumber
&&~t^T\nabla f_i(y_*)\le0, ~i\in\{j\in\{2,\cdots,k\}: f_j(y_*)=0, \mu^*_j=0\}\}.\nonumber
\end{eqnarray}
\item
Suppose $(x_*,y_*)$ is feasible to \eqref{eq:TRS-C} and there exist $\mu^*_j\ge 0$ ($j=1,\cdots,k$) satisfying \eqref{eq:Lx}, \eqref{eq:Ly}, \eqref{eq:Lc}, and
%\begin{equation*}
%\bmatrix
%s^T & t^T
%\endbmatrix
%\bmatrix
%H+\mu^*_1 I & 0 \\
%0 & \nabla^2 f_0(y_*)+\frac{\mu^*_1}{2} \nabla^2 f_1(y_*)+\sum_{j=2}^k\mu^*_j \nabla^2f_j(y_*)\\
%\endbmatrix
%\bmatrix
%s \\
%t \\
%\endbmatrix
%> 0,~\forall (s, t)\in S'\setminus\{0\},
%\end{equation*}
%holds for any $(s, t)$ satisfying $(s, t)\neq0, 2s^Tx_*+t^T\nabla f_1(y_*)=0$ and $t^T\nabla f_i(y_*)=0$ for $i=2,\cdots,n$ with $f_i(y_*)=0$.
\begin{eqnarray}
&&\bmatrix
s^T & t^T
\endbmatrix
\bmatrix
H+\mu^*_1 I & 0 \\
0 & \nabla^2 f_0(y_*)+\frac{\mu^*_1}{2} \nabla^2 f_1(y_*)+\sum_{j=2}^k\mu^*_j \nabla^2f_j(y_*)\\
\endbmatrix
\bmatrix
s \\
t \\
\endbmatrix
> 0,~\forall (s, t)\in T\setminus\{0\}.\nonumber
\end{eqnarray}
Then $(x_*,y_*)$ is a strict local minimizer of \eqref{eq:TRS-C}.
\end{enumerate}
\end{lemma}	
\emph{Proof.}
Firstly, we show that LICQ is true at $(x_*, y_*)$ under Assumption \ref{as:LICQ}.
If $x_*=0$, as $(0, \nabla f_i(y_*)^T)$ for $i\in\{ j\in\{ 1,\cdots,k\}:~ f_j(y_*)=0\}$ are linearly independent according to Assumption \ref{as:LICQ}, LICQ holds;  If $x_*\neq0$, $(0, \nabla f_i(y_*)^T)$ for $i\in\{ j\in\{ 2,\cdots,k\}:~ f_j(y_*)=0\}$ together with $(x_*^T, \nabla f_1(y_*)^T)$ are linearly independent according to Assumption \ref{as:LICQ} and the fact that $x_*\neq 0$.

With LICQ holding, the proof except \eqref{eq:qjjds} in \eqref{eq:Lc} follows from the classical optimality conditions in nonlinear programming (see \cite{Fletcher1987}, Chapter 9; \cite{Luenberger1984}, Chapter 11).
We conclude that \eqref{eq:qjjds} hold with the same reason as in Theorem \ref{theo:qjjbytj}.
\eproof
	
Next, we show that under Assumption \ref{as:LICQ}, for any
local non-global minimizer $(x_*,y_*)$, $x_*\neq 0$.

\begin{theorem}\label{th:local_peq0}
Suppose $(0,y_*)$ is a local minimizer of \eqref{eq:TRS-C}. Under Assumption \ref{as:LICQ},  $(0,y_*)$ is a global minimizer of \eqref{eq:TRS-C}.
\end{theorem}

\emph{Proof.}
According to Lemma \ref{le:traoptcon} (a), there exist $\mu^*_j\ge0$ ($j=1,\cdots,k$) such that \eqref{eq:Lx},\eqref{eq:Ly}, and \eqref{eq:Lc} hold. Furthermore, the inequality \eqref{eq:jbejbytj}
holds for any $(s,t)\in T$ in \eqref{st} with $x^*=0$.
Taking $t=0$ and $x^*=0$ into \eqref{eq:jbejbytj} and \eqref{st} gives
$$
s^T
(H+\mu^*_1 I)
s
\ge 0,~\forall s\in\R^n.
$$
Thus, $H+\mu^*_1 I\succeq 0$ and hence $(0, y_*)$ is a global minimizer of \eqref{eq:TRS-C} due to Theorem 2.
\eproof

Finally, we present additional necessary conditions of local non-global minimizer for \eqref{eq:TRS-C}. To this end, the following key lemma is required.

\begin{lemma}\label{le:M94}
\begin{itemize}
\item [(a)]
 {[Lemma 2.8, \cite{Sorensen82}]} If $x_*\in \R^n$ is a global minimizer of \eqref{eq:TRS}, there exists $\mu^*_1\ge \max\{-\lambda_1,0\}$ satisfying \eqref{eq:Lx}.
%(b) [Lemma 2.4,\cite{JOSE94}]
%If $x_*$ is a local minimizer of \eqref{eq:TRS}, then \eqref{eq:Lx} holds with $\mu^*_1\ge -\lambda_2$.

\item[(b)] {[Lemma 3.2,\cite{JOSE94}]}
If $c$ is orthogonal to some eigenvector associated with $\lambda_1$,  there is no local non-global minimizer of \eqref{eq:TRS}.

\item[(c)] {[Lemma 3.3,\cite{JOSE94}]}
If $x_*$ is a local non-global minimizer of \eqref{eq:TRS}, \eqref{eq:Lx} holds with $0\le\mu^*_1\in (-\lambda_2, -\lambda_1)$ and $x_*^Tx_*-\Delta=0$.

\item[(d)] {[Proposition 3.5, \cite{STEFANO98}]}
At the local non-global minimizer of \eqref{eq:TRS}, the strict complementarity condition holds.
\end{itemize}
\end{lemma}

\begin{theorem}\label{th:local_pneq0}
Let $(x_*,y_*)$ be a local non-global minimizer of \eqref{eq:TRS-C}. Under Assumption \ref{as:LICQ}, it holds that $\lambda_1<\lambda_2$,
there exist $\mu^*_j\ge0 $ ($j=1,\cdots,k$) such that  \eqref{eq:Lx}, \eqref{eq:Ly} and \eqref{eq:Lc} with $\max\{0,-\lambda_2\}<\mu^*_1<-\lambda_1$,
and $c$ is not orthogonal to any eigenvectors associated with $\lambda_1$.	
\end{theorem}

\emph{Proof.}
According to Lemma \ref{le:traoptcon} (a), there exist $\mu^*_j\ge0$  ($j=1,\cdots,k$) satisfying \eqref{eq:Lx}, \eqref{eq:Ly}, and \eqref{eq:Lc}.
Moreover, $x_*$ is a local minimizer of \eqref{eq:GTRSx}, since $(x_*,y_*)$ is locally optimal for \eqref{eq:TRS-C}. By Theorem \ref{th:local_peq0} and the fact that $(x_*,y_*)$ be a local non-global minimizer of \eqref{eq:TRS-C}, $x_*\neq0$. Then \eqref{eq:GTRSx} is an instance of \eqref{eq:TRS} with $\Delta=-f_1(y_*)>0$. We claim that $x_*$ is a local non-global minimizer of \eqref{eq:GTRSx}. Suppose this is not true, then $x_*$ is a global minimizer of \eqref{eq:GTRSx}. There exists $\bar\mu_1\ge -\lambda_1$ such that
\begin{equation}\label{eq:trcjbj1}
(H+\bar\mu_1 I)x_*+c=0
\end{equation}
according to Lemma \ref{le:M94} (a). Combining  \eqref{eq:Lx}, \eqref{eq:trcjbj1} with $x_*\neq 0$, we have $\mu^*_1=\bar\mu_1\ge -\lambda_1$.
Therefore, $(x_*,y_*)$ is globally optimal for \eqref{eq:TRS-C} due to Theorem \ref{theo:qjjcftj}. It contradicts the fact that $(x_*,y_*)$ is a local non-global minimizer of \eqref{eq:TRS-C}.

Since $x_*^T x_*+f_2( y_*)=0$ according to Lemma \ref{le:M94} (c), we get $\mu^*_1>0$ from Lemma \ref{le:M94} (d).
Then the remaining results to be proved follow from Lemma \ref{le:M94} (b), (c) and the fact that $x_*$ is a local non-global minimizer of \eqref{eq:GTRSx}.
\eproof

\section{Local Non-Global Optimality Conditions: Single-Constraint Case}\label{sec:5}

%To get more features for local non-global minimizers of \eqref{eq:TRS-C}, a simple case with $k=1$
In this section, we study local non-global optimality conditions of the single-constrained case of \eqref{eq:TRS-C}:
\begin{equation}
\begin{array}{cl}\label{eq:TRCS}
\min\limits_{x\in \R^n, y\in\R^m} & \tfrac{1}{2}x^THx+c^Tx+f_0(y)\\
\ST   & x^Tx+f(y)\le 0,
\end{array}
\end{equation}
where $f_0(y), f(y)$ are convex and twice continuously differentiable. Note that we do not make Assumption \ref{as:LICQ} in this section.

%Lemma \ref{lemma:trivial} will prove that it is trivial when $\{y: f(y)<0\}= \emptyset$.
We first identify a trivial case without local non-global minimizer.
\begin{remark}\label{lemma:trivial}
If $\{y: f(y)<0\}= \emptyset$ holds, \eqref{eq:TRCS} has no local nonglobal minimizer.
\end{remark}
\emph{Proof.}
If $(x_*, y_*)$ is a local minimizer of \eqref{eq:TRCS}, we have $f(y_*)=0$ and $x_*=0$, as $\{y: f(y)<0\}= \emptyset$ and $f(y_*)\le-x_*^Tx_*\le0$. Then
\eqref{eq:TRCS} reduces to the problem
$$
\min_{f(y)\le 0} f_0(y),
$$
which is a convex optimization problem. It follows that  \eqref{eq:TRCS} under the assumption $\{y: f(y)<0\}= \emptyset$ has no local non-global minimizer.
\eproof

%Next, we only consider the nontrivial case that
In the remainder of this section, without loss of generality, we assume
$\{y:f(y)<0\}\neq \emptyset$.
\begin{lemma}\label{LICQ-TRCS}
LICQ holds at any local nonglobal minimizer of \eqref{eq:TRCS}.
\end{lemma}
\emph{Proof.}
Let  $(x_*, y_*)$ be a local nonglobal minimizer of \eqref{eq:TRCS}. If $x_*\neq 0$, $(x_*^T, \nabla f(y_*)^T)\neq 0$, LICQ holds.
It's sufficient to consider the case $x_*=0$, which implies that $f(y_*)=0$ since $x_*^Tx_*+f(y_*)=0$. Under the assumption $\{y:f(y)<0\}\neq \emptyset$ and the convexity of $f(y)$, we have $\nabla f(y_*)\neq0$, which guarantees LICQ.
\eproof

Let $(x_*,y_*)$ be a local non-global minimizer of \eqref{eq:TRCS}.
According to Lemma \ref{LICQ-TRCS} and Theorem \ref{th:local_peq0}, $x_*\ne 0$. Moreover, it follows from Lemma \ref{LICQ-TRCS} and Theorem \ref{th:local_pneq0} that there is a unique $\mu^*\in (\max\{0,-\lambda_2\},-\lambda_1)$ such that
\begin{eqnarray}
&&(H+\mu^* I)x_*=-c,\label{eq:local_p}\\
&&\nabla f_0(y_*) + \frac{\mu^*}{2}\nabla f(y_*)= 0,\label{eq:local_w}\\
&&
x_*^T x_*+f(y_*)=0.\label{eq:local_phi}
\end{eqnarray}

Define the following scalar function similar to Mart\'{\i}nez \cite{JOSE94}
$$
\varphi(\mu) =\|{(H+\mu I)}^{-1}c\|^2,~\mu\in (-\lambda_2, -\lambda_1).
$$
Using \eqref{eq:tzzfjH},  one can verify that
\begin{equation}\label{eq:phi1ds}
\varphi(\mu)=\sum_{i=1}^n\frac{g_i^2}{{(\lambda_i+\mu)}^2},\
\varphi'(\mu)=-\sum_{i=1}^n\frac{2g_i^2}{{(\lambda_i+\mu)}^3},\ \text{and}\ \varphi''(\mu)=\sum_{i=1}^n\frac{6g_i^2}{{(\lambda_i+\mu)}^4}
\end{equation}
for all $\mu\in (-\lambda_2,-\lambda_1)$, where $g=V^Tc$. If a local non-global minimizer exists, we must have $g_1\ne 0$ and $\lambda_1<\lambda_2$ due to Theorem \ref{th:local_pneq0}.
%%%%%%%%%%
%\begin{assumption}\label{as:hsjzfqy}
%	$\nabla^2f_0(y_*)+\frac{\mu^*}{2}\nabla^2 f(y^*)$ is nonsingular.
%\end{assumption}

For easy analysis, we introduce the following nonsingular assumption, which automatically holds
if either  $f_0(y)$ or $f(y)$ is strongly convex.
\begin{assumption}\label{as:hsjzfqy}
At any $y\in\R^m$ with $f(y)<0$,  $\nabla^2f_0(y)+\frac{\mu_y}{2}\nabla^2 f(y)\succ 0$ for some $\mu_y>0$.
\end{assumption}
Assumption \ref{as:hsjzfqy} implies the following property.
\begin{lemma}\label{Ass2}
Suppose Assumption \ref{as:hsjzfqy} holds,
$\nabla^2f_0(y)+\frac{\mu}{2}\nabla^2 f(y)\succ 0$ for any $\mu>0$ and any $y$ such that $f(y)<0$.
\end{lemma}
\emph{Proof.}
Under Assumption \ref{as:hsjzfqy}, for any $y$ with $f(y)<0$ and any $\mu>0$, if $\mu\ge\mu_y$,
$$
\nabla^2f_0(y)+\frac{\mu}{2}\nabla^2f(y)\succeq \nabla^2f_0(y)+\frac{\mu_y}{2}\nabla^2f(y)\succ 0.
$$
Otherwise, $\mu<\mu_y$ and it holds that
$$\nabla^2f_0(y)+\frac{\mu}{2}\nabla^2 f(y)=\frac{\mu}{\mu_y}\left(\nabla^2f_0(y)+\frac{\mu_y}{2}\nabla^2f(y)\right)+\frac{\mu_y-\mu}{\mu_y}\nabla^2f_0(y)\succ0.
$$
The last inequality dues to the convexity of $f_0(y)$ and $\mu<\mu_y$.
The proof is complete.
\eproof

Let $(x_*, y_*)$ be a local non-global minimizer of \eqref{eq:TRCS}. We first have $\varphi(\mu^*)=-f(y_*)$ by \eqref{eq:local_p} and  \eqref{eq:local_phi}.
Since $x_*$ is a local non-global minimizer of \eqref{eq:GTRSx}, it holds that $\varphi'(\mu^*)\ge 0$ by Theorem 3.1 (i) of Mart\'{\i}nez \cite{JOSE94}. While in our case, we have a stronger result.

\begin{theorem}\label{th:local_nephi}
Suppose $(x_*,y_*)$ is a local non-global minimizer of \eqref{eq:TRCS}, there exist $\mu^*\in (\max\{0,-\lambda_2\}, -\lambda_1)$ such that \eqref{eq:local_p}, \eqref{eq:local_w} and \eqref{eq:local_phi} hold. Furthermore, under Assumption \ref{as:hsjzfqy}, we have
\begin{equation}\label{eq:phids}
\varphi'(\mu^*)\ge \frac{1}{2}\nabla f(y_*)^T\left[\nabla^2 f_0(y_*)+\frac{\mu^*}{2}\nabla^2 f(y_*)\right]^{-1}\nabla f(y_*).
\end{equation}
\end{theorem}

\emph{Proof.}
Theorem \ref{th:local_peq0} implies $x_*\neq 0$. According to Theorem \ref{th:local_pneq0}, there is a unique
$\mu^*$ such that \eqref{eq:local_p}, \eqref{eq:local_w} and \eqref{eq:local_phi} hold with $\mu^*\in (\max\{0,-\lambda_2\}, -\lambda_1)$.
Moreover, According to Lemma \ref{le:traoptcon} (a), we have
\begin{equation}\label{eq:ejbytj}
\bmatrix
s^T & t^T
\endbmatrix
\bmatrix
H+\mu^* I & 0\\
0 & \nabla^2 f_0(y_*)+\frac{\mu^*}{2}\nabla^2 f(y_*)
\endbmatrix
\bmatrix
s \\ t
\endbmatrix
\ge 0
\end{equation}
for all $s\in \R^n, t\in \R^m$ satisfying
\begin{equation}\label{eq:xxhkxfx}
x_*^Ts+\tfrac{1}{2}\nabla f(y_*)^Tt=0.
\end{equation}

Since $H+\mu^* I$ is nonsingular, it follows from \eqref{eq:tzzfjH} and \eqref{eq:local_p} that
\begin{equation}\label{eq:local_x}
x_*=-(H+\mu^* I)^{-1}c=-\sum_{i=1}^n\frac{g_i}{\lambda_i+\mu^*}v_i.
\end{equation}
For any $\mu\in (-\lambda_2,-\lambda_1)$,
define $W(\mu)\in\R^{(n+m)\times(n+m-1)}$:
\begin{equation}\label{w(u)}
W(\mu)=\bmatrix V & 0\\ 0 & I_m\endbmatrix
\bmatrix
    \frac{g_2}{\lambda_2+\mu} & \frac{g_3}{\lambda_3+\mu} & \cdots & \frac{g_n}{\lambda_n+\mu}
     & -\frac{1}{2}\frac{\partial f(y_*)}{\partial y_1}  & -\frac{1}{2}\frac{\partial f(y_*)}{\partial y_2} & \cdots & -\frac{1}{2}\frac{\partial f(y_*)}{\partial y_m}\\

     -\frac{g_1}{\lambda_1+\mu} & 0 & \cdots & 0
     & 0 & 0 & \cdots & 0\\

     0 & -\frac{g_1}{\lambda_1+\mu} & \cdots & 0
     & 0 & 0 & \cdots & 0\\

     \vdots & \vdots &   & \vdots
     & \vdots & \vdots &   & \vdots\\

     0 & 0 & \cdots & -\frac{g_1}{\lambda_1+\mu}
     & 0 & 0 & \cdots & 0\\

     0 & 0 & \cdots & 0
     & -\frac{g_1}{\lambda_1+\mu} & 0 & \cdots & 0\\

     0 & 0 & \cdots & 0
     & 0 & -\frac{g_1}{\lambda_1+\mu} & \cdots & 0\\

     \vdots & \vdots &   & \vdots
     & \vdots & \vdots &   & \vdots\\

     0 & 0 & \cdots & 0
     & 0 & 0 & \cdots & -\frac{g_1}{\lambda_1+\mu}\\
\endbmatrix.
\end{equation}
By Theorem \ref{th:local_pneq0}, we have $g_1=v_1^Tc\ne 0$. Therfore,
\begin{equation}\label{eq:xxkxfxR}
{\rm rank} (W(\mu))=n+m-1,~\forall \mu\in (-\lambda_2,-\lambda_1).
\end{equation}
Moreover, by \eqref{eq:local_x} and \eqref{w(u)} we have
\begin{equation}\label{eq:xxhkxfxj}
  \bmatrix   x_*^T, &  \tfrac{1}{2}\nabla f(y_*)^T   \endbmatrix
  W(\mu^*)e_j=0,~j=1,\cdots,n+m-1,
\end{equation}
where $e_j$ is the $j$-th standard unit vector in $\R^{n+m-1}$.
According to \eqref{eq:xxkxfxR} and \eqref{eq:xxhkxfxj}, all columns of $W(\mu^*)$ form a basis of the hyperplane \eqref{eq:xxhkxfx}.
For any $\mu\in (-\lambda_2,-\lambda_1)$, define
\begin{equation}\label{eq:jyhsjz}
B(\mu)=W(\mu)^T
\bmatrix   H+\mu I & 0\\
          0 & \nabla^2 f_0(y_*)+\frac{\mu^*}{2}\nabla^2 f(y_*)
          \endbmatrix
W(\mu).
\end{equation}
Then, we have
\begin{equation*}\label{eq:jyhsjz1}
B(\mu)=
\bmatrix   \hat B(\mu)+\sigma uu^T & -\tfrac{\sigma}{2}u\nabla f(y_*)^T\\
          -\tfrac{\sigma}{2}\nabla f(y_*)u^T & ~ \frac{g_1^2}{\sigma^2}\left(\nabla^2 f_0(y_*)+\frac{\mu^*}{2}\nabla^2 f(y_*)\right)+\frac{\sigma}{4}\nabla f(y_*)\nabla f(y_*)^T
          \endbmatrix
\end{equation*}
where
\begin{equation*}
\hat B(\mu)=\bmatrix
\frac{g_1^2(\lambda_2+\mu)}{{(\lambda_1+\mu)}^2}& \cdots &0\\
\vdots               & \ddots  & \vdots\\
0      & \cdots   &\frac{g_1^2(\lambda_n+\mu)}{{(\lambda_1+\mu)}^2}
\endbmatrix
,~
u=\bmatrix
\frac{{g_2}}{\lambda_2+\mu}\\
\vdots\\
\frac{{g_n}}{\lambda_n+\mu}\\
\endbmatrix
,~ \sigma=\lambda_1+\mu.
\end{equation*}
Equivalently,
\begin{equation}\label{eq:jbfqj6}
B(\mu) =\bar B(\mu) +\sigma \bar u\bar u^T,
\end{equation}
where
$$
\bar B(\mu)=\bmatrix
  \hat B(\mu) & 0\\
  0 & (\frac{g_1}{\sigma})^2(\nabla^2 f_0(y_*)+\frac{\mu^*}{2}\nabla^2 f(y_*))
\endbmatrix,~
\bar u=\bmatrix
  u\\
 -\frac{1}{2}\nabla f(y_*)
\endbmatrix.
$$
As $g_1\neq 0$, $\hat B(\mu)\succ0$ for all $\mu\in (-\lambda_2, -\lambda_1)$.  Since both $f_0(\cdot)$ and $f(\cdot)$ are convex and $\mu^*> 0$, according to Assumption \ref{as:hsjzfqy} and Lemma \ref{Ass2}, we have
\begin{equation*}
\nabla^2 f_0(y_*)+\frac{\mu^*}{2}\nabla^2 f(y_*)\succ 0.
\end{equation*}
Now, $\bar B(\mu)\succ0$ for all $\mu\in (-\lambda_2.-\lambda_1) $. So
\begin{align}\label{eq:jbfqj2}
&{\rm det}(\bar B(\mu))>0.
\end{align}
Hence, $B(\mu)=\bar B(\mu)(I+\sigma\bar B(\mu)^{-1}\bar u\bar u^T)$ and (see \cite{GOLUB89})
\begin{equation}\label{eq:jbfqj1}
{\rm det}(B(\mu))={\rm det}(\bar B(\mu))(1+\sigma\bar u^T[\bar B(\mu)]^{-1}\bar u).
\end{equation}
Notice that
\begin{align}
&1+\sigma \bar u^T[\bar B(\mu)]^{-1}\bar u \nonumber\\
 =& 1 + \sigma u^T\hat B^{-1}(\mu)u + \frac{\sigma^3}{4g_1^2}\nabla f(y_*)^T\left[\nabla^2 f_0(y_*)+\frac{\mu^*}{2}\nabla^2 f(y_*)\right]^{-1}\nabla f(y_*) \nonumber\\
=&-\frac{\sigma^3}{2g_1^2}\left[\frac{-2g_1^2}{(\lambda_1+\mu)^3}+\cdots+\frac{-2g_n^2}{(\lambda_n+\mu)^3}- \frac{1}{2}\nabla f(y_*)^T\left[\nabla^2 f_0(y_*)+\frac{\mu^*}{2}\nabla^2 f(y_*)\right]^{-1}\nabla f(y_*)\right] \nonumber\\[0.2cm]
=&-\frac{\sigma^3}{2g_1^2}\left[\varphi'(\mu)- \frac{1}{2}\nabla f(y_*)^T\left[\nabla^2 f_0(y_*)+\frac{\mu^*}{2}\nabla^2 f(y_*)\right]^{-1}\nabla f(y_*)\right], \label{eq:jbfqj3}
\end{align}
where the third equality follows from \eqref{eq:phi1ds}. It follows from \eqref{eq:jbfqj1} and \eqref{eq:jbfqj3} that
\begin{equation}\label{eq:jbfqj4}
{\rm det}(B(\mu))=-\frac{\sigma^3}{2g_1^2}{\rm det}(\bar B(\mu))\left[\varphi'(\mu)- \frac{1}{2}\nabla f(y_*)^T\left[\nabla^2 f_0(y_*)+\frac{\mu^*}{2}\nabla^2 f(y_*)\right]^{-1}\nabla f(y_*)\right].
\end{equation}
According to the definition of $\sigma$ and the fact $g_1\neq 0$, for $\mu^*\in (-\lambda_2,-\lambda_1)$, it holds that
\begin{equation}\label{eq:jbfqj5}
-\frac{\sigma^3}{2g_1^2}>0.
\end{equation}
It implies from \eqref{eq:ejbytj}, \eqref{eq:xxhkxfx} and \eqref{eq:jyhsjz} that $B(\mu^*)\succeq 0$ and hence ${\rm det}(B(\mu^*))\ge 0$. Consequently, by \eqref{eq:jbfqj2}, \eqref{eq:jbfqj4} and \eqref{eq:jbfqj5}, we get \eqref{eq:phids}.
\eproof

Based on the above defined $\varphi(\mu)$, we now  characterize the second-order sufficient optimality condition.

\begin{theorem}\label{th:local_suphi}
Under Assumption \ref{as:hsjzfqy},   if there is a $\mu^*\in (\max\{0,-\lambda_2\}, -\lambda_1)$ satisfying
\eqref{eq:local_p}-\eqref{eq:local_phi},
 and
\begin{equation}\label{eq:phidsdy}
\varphi'(\mu^*)>\frac{1}{2}\nabla f(y_*)^T\left[\nabla^2 f_0(y_*)+\frac{\mu^*}{2}\nabla^2 f(y_*)\right]^{-1}\nabla f(y_*),
\end{equation}
then  $(x_*,y_*)$ is a strict local minimizer of \eqref{eq:TRCS}.
\end{theorem}

\emph{Proof.}
%We first prove that $B(\mu^*)\succ0$.
Assume that $(x_*,y_*)$ satisfies \eqref{eq:local_p}-\eqref{eq:local_phi} for some $\mu^*\in (\max\{0,-\lambda_2\}, -\lambda_1)$ and \eqref{eq:phidsdy} holds. Since $\varphi'(\cdot)$ is strictly increasing in $(-\lambda_2, -\lambda_1)$, we have
\begin{equation}\label{eq:jbfqjcf1}
\varphi'(\mu)> \frac{1}{2}\nabla f(y_*)^T\left[\nabla^2 f_0(y_*)+\frac{\mu^*}{2}\nabla^2 f(y_*)\right]^{-1}\nabla f(y_*),~\forall \mu\in [\mu^*, -\lambda_1).
\end{equation}

Suppose that $(x_*,y_*)$ is not a strictly local minimizer of \eqref{eq:TRCS}.
Define $W(\mu)$ and $B(\mu)$ as in \eqref{w(u)} and \eqref{eq:jyhsjz}, respectively. Notice that \eqref{eq:jbfqj4} and \eqref{eq:jbfqj5} hold for all $\mu\in [\mu^*, -\lambda_1)$.  All the eigenvalues of $\bar B(\mu)$ are positive for $\mu\in (-\lambda_2,-\lambda_1)$.
If $(x_*,y_*)$ is not a strict local minimizer of \eqref{eq:TRCS}, according to Lemma \ref{le:traoptcon} (b), $B(\mu^*)$ has at least one eigenvalue less than or equal to zero. Now by \eqref{eq:jbfqj6}, all the eigenvalues of $B(\mu)$ are strictly positive if $\mu\in (-\lambda_2,-\lambda_1)$ is close enough to $-\lambda_1$. Therefore, there exists $\tilde\mu\in[\mu^*,-\lambda_1)$ such that $B(\tilde\mu)$ is singular. So,
\begin{equation}\label{eq:jbfqjcf2}
0={\rm det}(B(\tilde\mu))=-\frac{\tilde\sigma^3}{2g_1^2}{\rm det}(\bar B(\tilde\mu))\left[\varphi'(\tilde\mu)- \frac{1}{2}\nabla f(y_*)^T\left[\nabla^2 f_0(y_*)+\frac{\mu^*}{2}\nabla^2 f(y_*)\right]^{-1}\nabla f(y_*)\right],
\end{equation}
where $\tilde\sigma=\lambda_1+\tilde\mu$. It implies from \eqref{eq:jbfqjcf2} and ${\rm det}(\bar B(\tilde\mu))>0$ that
$$
\varphi'(\tilde\mu)=\frac{1}{2}\nabla f(y_*)^T\left[\nabla^2 f_0(y_*)+\frac{\mu^*}{2}\nabla^2 f(y_*)\right]^{-1}\nabla f(y_*),
$$
which
contradicts \eqref{eq:jbfqjcf1}. Thus, the fact $\mu^*\in (\max\{0,-\lambda_2\}, -\lambda_1)$,  \eqref{eq:local_p}-\eqref{eq:local_phi}, and \eqref{eq:phidsdy} imply that $(x_*, y_*)$ is a strict minimizer of \eqref{eq:TRCS}.
\eproof

%\begin{assumption}\label{as:hsjzqt}
%	Either $f_0$ or $f$ is \textcolor{red}{strong} convex.
%\end{assumption}

Under Assumption \ref{as:hsjzfqy}, we show how to find the local non-global minimizer of \eqref{eq:TRCS}. For $\mu\in (\max\{-\lambda_2, 0\},-\lambda_1)$, define
\begin{equation*}
x(\mu)=-(H+\mu I)^{-1}c,~ y(\mu)=\arg\min\limits_{y\in\R^m} \ f_0(y)+\frac{\mu}{2}f(y).
\end{equation*}
Then, we have
\begin{equation}\label{eq:jbfqjy22}
\nabla f_0(y(\mu)) = -\frac{\mu}{2}\nabla f(y(\mu)).
\end{equation}
Notice that, under Assumption \ref{as:hsjzfqy}, $y(\mu)$  is uniquely defined and continuously differential  in terms of $\mu$, according to the well-known implicit function theorem. Denote by $y'(\mu)$ the derivative of $y(\mu)$.
Differentiating both sides of equation \eqref{eq:jbfqjy22} with respect to $\mu$  yields that
$$
\left[\nabla^2 f_0(y(\mu))+\frac{\mu}{2}\nabla^2 f(y(\mu))\right]y'(\mu)=-\frac{1}{2}\nabla f(y(\mu)).
$$
Note that Assumption \ref{as:hsjzfqy} and Lemma \ref{Ass2} imply the non-singularity of $\nabla^2 f_0(y(\mu))+\frac{\mu}{2}\nabla^2 f(y(\mu))$. Then we have
\begin{equation}\label{eq:gw(mu)}
y'(\mu)=-\frac{1}{2}{[\nabla^2 f_0(y(\mu))+\frac{\mu}{2}\nabla^2 f(y(\mu))]}^{-1}\nabla f(y(\mu)).
\end{equation}
Define
\begin{equation*}
\psi(\mu)=-f(y(\mu)).
\end{equation*}
Applying the chain rule to calculate the derivative of $\psi(\cdot)$ gives
$$
\psi'(\mu)= -\nabla f(y(\mu))^Ty'(\mu)= \frac{1}{2}\nabla f(y(\mu))^T{[\nabla^2 f_0(y(\mu))+\frac{\mu}{2}\nabla^2 f(y(\mu))]}^{-1}\nabla f(y(\mu)),
$$
where the second equality holds due to \eqref{eq:gw(mu)}. According to Theorem \ref{th:local_suphi}, if there exists $\mu\in (\min\{-\lambda_2,0\},-\lambda_1)$ such that
$$
\varphi(\mu)=\psi(\mu), ~ \varphi'(\mu)>\psi'(\mu),
$$
then $(x(\mu), y(\mu))$ is a local non-global minimizer of \eqref{eq:TRCS}. It is sufficient to find the root of  the scalar function
\begin{equation}\label{eq:scalarF}
\phi(\mu)=\varphi(\mu)-\psi(\mu),~ \mu\in (\max\{-\lambda_2, 0\},-\lambda_1)
\end{equation}
satisfying $\phi'(\mu)>0$.

%%%%%%%%%%%%%%%%%%%%%%%%%%%%%%%%%%%%%%%%%%%%%5
%\section{\eqref{eq:TRS-C} could have more than one local non-global minimizer.}\label{sec:6}
\section{Local Non-Global Optimality Conditions: an Intensive Analysis on Quadratic Single-Constraint Case.}\label{sec:6}
In this section, we focus on a more special single-constraint case of \eqref{eq:TRS-C}:
\begin{equation}\label{eq:TRS-L}
\begin{array}{cl}
\min\limits_{x\in \R^n, y\in\R} & \tfrac{1}{2}x^THx+c^Tx+f_0(y)\\
\ST   & x^Tx-ay-b\le 0,
\end{array}
\end{equation}
where $f_0$ is twice continuously differentiable and strongly convex in $S:=\{y:ay+b> 0\}$. We show that in this case, there may be more than one local non-global minimizer, and under some assumptions there is a
necessary and sufficient optimality condition for local non-global minimizer.

Specially structured, \eqref{eq:TRS-L} contains two well-known cases, \eqref{eq:TRS} and \eqref{eq:pRS-gtrs}.
Specifically, \eqref{eq:TRS-L} with $a=0, b=\Delta, f_0(y)=y^2$ reduces to \eqref{eq:TRS}, and  \eqref{eq:pRS-gtrs} corresponds to the case of \eqref{eq:TRS-L}  with $a=1, b=0, f_0(y)=\tfrac{\sigma}{p}y^{\tfrac{p}{2}}$. A common property of both cases is that there is at most one local non-global minimizer, characterized by a necessary and sufficient condition \cite{JOSE94,Xia17,Jiulin20}.

Without loss of generality, we can assume that $a>0$.
Actually, if $a=0$,  \eqref{eq:TRS-L} separates into (TRS) and an unconstrained convex optimization problem in terms of $x$ and $y$, respectively.
If $a<0$. let $z=-y$. Then \eqref{eq:TRS-L} is equivalent to
\begin{equation*}
\begin{array}{cl}
\min\limits_{x\in \R^n, z\in\R} & \tfrac{1}{2}x^THx+c^Tx+f_0(-z)\\
\ST   & x^Tx-(-a)z-b\le 0,
\end{array}
\end{equation*}
where $f_0(-z)$ is convex in terms of $z$.

%\begin{assumption}\label{as:jbfqjjgs}
%	For given $a>0, b\ge 0$, $f_0$ is twice continuous differentiable and \textcolor{red}{strong} convex in $S:=\{y | ay+b> 0\}$.
%\end{assumption}

As shown in Section 2, there are necessary and sufficient optimality conditions for  global minimizers of \eqref{eq:TRS-L} (see Theorems \ref{theo:qjjbytj} and \ref{theo:qjjcftj}). According to Theorem \ref{th:local_nephi}, for any local non-global minimizer of \eqref{eq:TRS-L} denoted by $(x_*,y_*)$, we have $x_*\neq 0$, and there is a unique $\max\{0,-\lambda_2\}<\mu^*<-\lambda_1$ such that
\begin{equation*}\label{eq:trsl_p}
(H+\mu^* I) x_*=-c,~
f_0'(y_*) - \frac{\mu^*a}{2} = 0,~
x_*^T x_*-ay_*-b=0,~\varphi'(\mu^*)\ge \frac{a^2}{2f''_0(y_*)}.
\end{equation*}
Since $f_0(\cdot)$ is strongly convex and twice continuously differentiable in $S$, we have $f''_0(\cdot)>0$, and $f'_0(\cdot)$ is strictly increasing so that  $(f'_0)^{-1}(\frac{\mu a}{2})$ exists.
Define
\begin{equation}\label{eq:zjbfqj-y}
y(\mu)={(f'_0)}^{-1}\left(\frac{\mu a}{2}\right).
\end{equation}
Finding a local non-global minimizer of \eqref{eq:TRS-L} amounts to search in $(\max\{-\lambda_2,0\},-\lambda_1)$ a root of $\varphi(\mu)=\psi(\mu)$  such that $\varphi'(\mu)>\psi'(\mu)$, where
\begin{equation}\label{eq:zjbfqj-psi}
\psi(\mu)=ay(\mu)+b=a{(f'_0)}^{-1}\left(\frac{\mu a}{2}\right)+b.
\end{equation}
Then, by Theorem \ref{th:local_suphi}, $(x(x),y(\mu))$ is a local non-global minimizer of \eqref{eq:TRS-L}.

%We care about whether there are necessary and sufficient optimality conditions for \eqref{eq:TRS-L} and whether there is at most one local non-global minimizer like \eqref{eq:TRS} and \eqref{eq:pRS-gtrs}.

While either \eqref{eq:TRS} or \eqref{eq:pRS-gtrs} has at most one local non-global minimizer, the answer to
Question \ref{ques2} could be surprisingly negative.
The following two examples illustrate that \eqref{eq:TRS-L} could have more than one local non-global minimizer.

\begin{example}[A quartic example of \eqref{eq:TRS-L} with two local non-global minimizers]\label{ex:l1}
Consider \eqref{eq:TRS-L} with $n=2, H={\rm Diag}\{-5,-1\}, c=(1,1), a=1,b=0$ and $f_0(y)$ is a quartic polynomial function:
\begin{equation*}
f_0(y)=\left(\tfrac{12377}{51072} - \tfrac{25\sqrt{210}}{3648}\right)y^4 + \left(\tfrac{5\sqrt{210}}{228} - \tfrac{9257}{7980}\right)y^3 + \left(\tfrac{1366171}{638400} - \tfrac{35\sqrt{210}}{1824}\right)y^2
+\left(\tfrac{\sqrt{210}}{190} + \tfrac{4667}{26600}\right)y.
\end{equation*}
One can verify the strong convexity of $f_0$. It follows from  \eqref{eq:phi1ds} and \eqref{eq:zjbfqj-psi} that
\begin{eqnarray}
&&\varphi(\mu)=\frac{1}{{(\mu-5)}^2}+\frac{1}{{(\mu-1)}^2},
\label{eq:zjbfqj-phil1}\\
&&\psi(\mu)=y(\mu)={(f'_0)}^{-1}\left(\tfrac{\mu}{2}\right).
\nonumber
\end{eqnarray}
We plot in Figure \ref{phi1Andphi1} the functions $\varphi(\mu)$ and $\psi(\mu)$. It is observed that, in the interval $(1,5)$, there are four roots of $\varphi(\mu)= \psi(\mu)$, $\mu_1=3.13,\mu_2=3.72,\mu_3=4.17,\mu_4=4.25$. As $\varphi'(\mu_i)>\psi'(\mu_i)$ for $i=2,4$,  $(x(\mu_2)^T, y(\mu_2))$ and $(x(\mu_4)^T, y(\mu_4))$  are two local non-global minimizers of Example \ref{ex:l1} by Theorem \ref{th:local_suphi}.

Actually, at $\mu_2$ and $\mu_4$, the corresponding solutions are $(x(\mu_2)^T, y(\mu_2))=(0.78, -0.37, 0.74)$ and $(x(\mu_4)^T, y(\mu_4))=(1.34, -0.31, 1.89)$, respectively. The reduced Hessian matrices \eqref{eq:jyhsjz} at $\mu_2, \mu_4$ are given by
\begin{equation*}
B(\mu_2)=\bmatrix
1.48 & -0.24\\
-0.24  & 0.24\\
\endbmatrix\succ 0,~
B(\mu_4)=\bmatrix
5.77 & -0.11\\
-0.11  & 0.36\\
\endbmatrix\succ 0.
\end{equation*}
Therefore, according to Lemma \ref{le:traoptcon} (b), both $(x(\mu_2)^T, y(\mu_2))$ and $(x(\mu_4)^T, y(\mu_4))$ are local minimizers of Example \ref{ex:l1}. On the other hand, there is a root $\mu_0=5.63$ of $\varphi(\mu)= \psi(\mu)$ in the interval $(5,+\infty)$. The corresponding solution $(x(\mu_0)^T, y(\mu_0))=(-1.58, -0.22, 2.56)$ is a global minimizer by Theorem \ref{theo:qjjcftj}.
\begin{figure}[h]
\centering
\includegraphics[width=15cm,height=7cm]{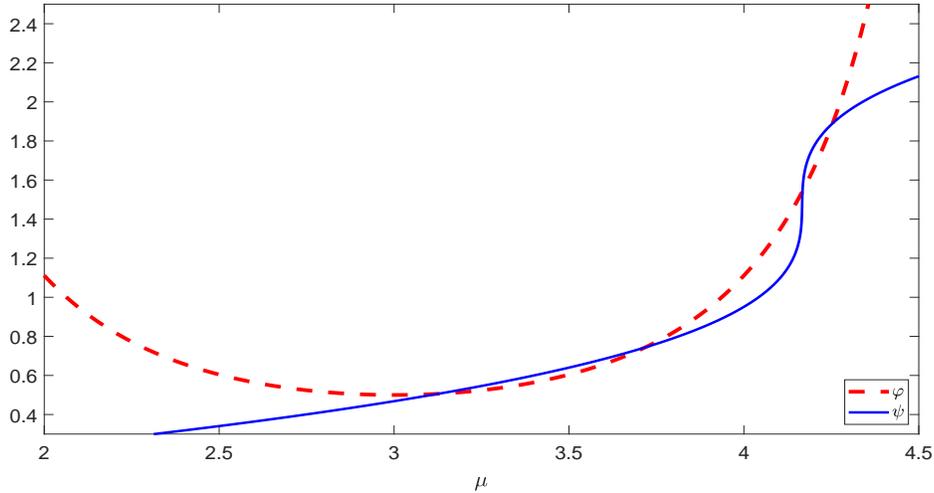}
\caption{Variations of $\varphi(\mu)$ and $\psi(\mu)$ in Example \ref{ex:l1}.}
\label{phi1Andphi1}
\end{figure}
\end{example}

In the following, we show that \eqref{eq:TRS-L} could have arbitrary number of local non-global minimizers.

\begin{example}[An example of \eqref{eq:TRS-L} with $d$ local non-global minimizers]\label{ex:l2}
Let $n=2, H={\rm Diag}\{-5,-1\}, c=(1,1), a=1,b=0$, and $d$ is a given positive integer.
%a $f_0(y)$ is built such that the resulting example has $d$ local non-global minimizers.
We first select a strictly increasing sequence $\{\mu_1,\cdots,\mu_{2d}\}$ from $(\max\{0,-\lambda_2\},-\lambda_1)$. We assume that $\mu_1$ is sufficiently close to $-\lambda_1$ so that $\varphi(\mu)$
strictly increases when $\mu>\mu_1$.
For each $j\in\{1,\cdots,d\}$, let $L_j(\mu):~\R\rightarrow \R$ be a linear function passing through two points  $(\mu_{2j-1},\varphi(\mu_{2j-1}))$ and $(\mu_{2j},\varphi(\mu_{2j}))$. Define
$$
L(\mu) = \max\{L_1(\mu),\cdots,L_d(\mu)\},
$$
which is piecewise linear and  convex in terms of $\mu$.  %We smooth $L(\mu)$ at the dis
$L(\mu)$ is nonsmooth at each intersection point of $L_j(\mu)$ and $L_{j+1}(\mu)$, denoted by $o_j$, for $j=1,\cdots,d-1$. Let
% is smoothed by replacing $L$ by quadratic functions in the specific small intervals.
%Denote by $o_j$ the $\mu$-coordinate of the intersection of the plots of $L_j$ and $L_{j+1}$. Let
$$
\epsilon=\tfrac{1}{2}
\min\{\min\{|o_j-\mu_{2j}|,|o_j-\mu_{2j+1}|\}:~ j=1,\cdots,d-1\}.
$$
Define $l_j=o_j-\epsilon, ~r_j=o_j+\epsilon$. For $j=1,\cdots,d-1$,
let $Q_j(\mu):~\R\rightarrow \R$  be a quadratic function not only connecting two endpoints $(l_j,L_{j}(l_j))$ and $(r_j,L_{j+1}(r_j))$,
but also tangent to $L_{j}(\mu)$ and $L_{j+1}(\mu)$ at these two endpoints, respectively. Define
$$
\psi(\mu)=\left\{\begin{array}{lll}
L_j(\mu), & \mu\in [r_{j-1},l_j],&\text{for}~j=1,\cdots, d,\\
Q_j(\mu), & \mu\in [l_j,r_j],&\text{for}~j=1,\cdots, d-1,
\end{array}
\right.
$$
where $r_0:=-\lambda_2$ and $l_d:=-\lambda_1$. We can verify that $\psi(\mu)$ is continuously differentiable and strictly increasing.

With the above definitions, we can see that  $\psi(\mu)$ intersects $\varphi(\mu)$ at $(\mu_i,\varphi(\mu_i)), i=1,2,\cdots,2d$, and it holds that $\varphi'(\mu_{2j})>\psi'(\mu_{2j})=L'_j(\mu_{2j})$ for $j=1,\cdots,d$.
According to \eqref{eq:zjbfqj-psi} and $(a,b)=(1,0)$, we have $y(\mu)=\psi(\mu)$. It follows from  $f_0'(y)=\tfrac{\mu}{2}$ that
\begin{equation*}
f_0'(y)=\frac{1}{2}\psi^{-1}(y),\label{eq:mbyds}
\end{equation*}
which implies that
\begin{equation*}
f_0(y)=\frac{1}{2}\int \psi^{-1}(y)dy
\end{equation*}
is convex, and $f''_0(y)>0$ in $(-\lambda_2,-\lambda_1)$  as $\psi^{-1}(y)$ is strictly increasing.
Since $\varphi'(\mu_{2j})>\psi'(\mu_{2j})$ for $j=1,\cdots,d$,
by Theorem \ref{th:local_suphi}, $(x(\mu_{2j}),\psi(\mu_{2j}))$ ($j=1,\cdots d$) are all local non-global minimizers.

%An example is presented to illustrate the general building method with $d=3$, which has three local non-global minimizers.

When $d=3$, setting $\mu_1=3.00,\mu_2=3.58,\mu_3=3.94,\mu_4=4.13,\mu_5=4.40,\mu_6=4.45,o_1=3.80,o_2=4.30$
in the above scheme gives
\begin{equation*}
	f_0(y)=\left\{
	\begin{aligned}
	y^2+\frac{y}{2},  & &y\in\left(-\infty,\frac{27}{40}\right] \\ \frac{1433y}{780}+\frac{2\cdot10^{1/2}{{(195y-131)}^{3/2}}}{114075}-\frac{3263447}{7300800}, & &y\in\left(\frac{27}{40},\frac{23}{25}\right]\\
	\frac{5y^2}{44}+\frac{383y}{200}-\frac{38363}{88000},  & &y\in\left(\frac{23}{25},\frac{79}{50}\right]\\
	\frac{871y}{420}+ \frac{{(2100y-3197)}^{3/2}}{1323000}-\frac{7189183}{10584000}, & &y\in\left(\frac{79}{50},\frac{143}{50}\right]\\ \frac{5y^2}{212}+\frac{2189y}{1060}-\frac{206802453584843}{281474976710656},  & &y\in\left(\frac{143}{50},\infty\right)\\
	\end{aligned}
	\right..
	\end{equation*}
In this case, $\varphi(\mu)$ is given in \eqref{eq:zjbfqj-phil1} and
\begin{equation*}
	\psi(\mu)=y(\mu)=\left\{
	\begin{aligned}
	\frac{\mu}{4}-\frac{1}{4},  & & \mu\in\left(-\infty,\frac{37}{10}\right] \\
	\frac{39{\mu}^2}{8}-\frac{1433}{40}\mu+\frac{53191}{800}, & & \mu\in\left(\frac{37}{10},\frac{39}{10}\right]\\
	\frac{11}{5}\mu-\frac{383}{50},  & &\mu\in\left(\frac{39}{10},\frac{42}{10}\right) \\
	{{21}{\mu}^2}-\frac{871}{5}\mu+\frac{18139}{50}, & &
	\mu\in\left(\frac{42}{10},\frac{44}{10}\right]\\
	\frac{53}{5}\mu-\frac{2189}{50},  & & \mu\in\left(\frac{44}{10},\infty\right) \\
	\end{aligned}
	\right..
\end{equation*}
We plot both in Figure \ref{phi1Andphi2}. We observe that $\varphi(\mu_j)= \psi(\mu_j)$ for $j=1,2,\cdots,6$ and $\varphi'(\mu_{2j})>\psi'(\mu_{2j})$ for $j=1,2,3$. Thus,  $(x(\mu_{2j}), \psi(\mu_{2j}))$ ($j=1,2, 3$) are three local non-global minimizers of Example \ref{ex:l2} according to Theorem \ref{th:local_suphi}.
\begin{figure}[h]
		\centering		\includegraphics[width=16cm,height=8cm]{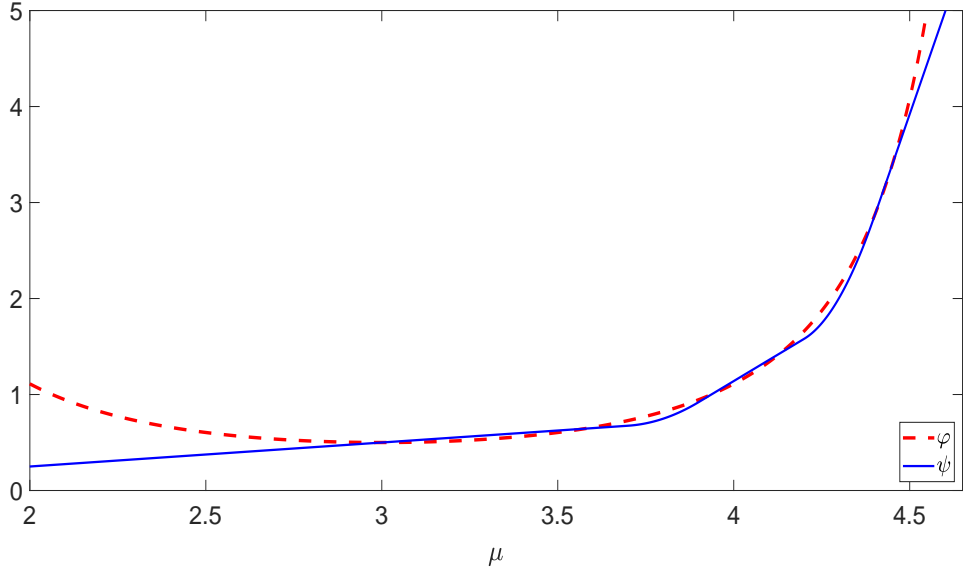}
		\caption{Variations of $\varphi(\mu)$ and $\psi(\mu)$ in  Example \ref{ex:l2}.}
		\label{phi1Andphi2}
	\end{figure}	
\end{example}
%Thus, by \eqref{eq:mbyds},
%	\begin{equation*}
%	f'_0(y)=\frac{\psi^{-1}(y)}{2}=\left\{
%	\begin{aligned}
%	2y+\frac{1}{2},  & & y\in\left(-\infty,\frac{27}{40}\right] \\
%	\frac{1433}{780}+\frac{2{{(39y/2-131/10)}^{1/2}}}{39}, & & y\in\left(\frac{27}{40},\frac{23}{25}\right]\\
%	\frac{5y}{22}+\frac{383}{220},  & & y\in\left(\frac{23}{25},\frac{79}{50}\right] \\
%	\frac{871}{420}+\frac{{{(84y-3197/25)}^{1/2}}}{84}, & & y\in\left(\frac{79}{50},\frac{143}{50}\right]\\
%	\frac{5y}{106}-\frac{2189}{1060},  & &y\in\left(\frac{143}{50},\infty\right)\\
%	\end{aligned}
%	\right.
%	\end{equation*}
%Furthermore, by \eqref{eq:mby}, we have the twice continuously differential and convex

%Examples \ref{ex:l1} and  \ref{ex:l2} give a negative answer to Question \ref{ques2}.

Now, we are interested in when \eqref{eq:TRS-L} has at most one local non-global minimizer, and whether\eqref{eq:TRS-L} could have necessary and sufficient optimality condition at local non-global minimizer.
The proof is a combination of those of Theorem 3.2 in \cite{Xia17} and  Theorem 3.1 in \cite{Jiulin20}.
%and Theorem \ref{local_sune} consider whether exists the sufficient and necessary conditions for local non-global minimizer of \eqref{eq:TRS-L}.

\begin{theorem}\label{local_onlyone}
Suppose $f_0(y)$ is thrice continuously differentiable and strongly convex in $S$ and $\psi(\mu)$ defined in \eqref{eq:zjbfqj-psi} is log-concave in $(\max\{-\lambda_2,0\},-\lambda_1)$, then

(i) \eqref{eq:TRS-L} has at most one local non-global minimizer.

(ii) $(x_*,y_*)$ is a local non-global minimizer of \eqref{eq:TRS-L} if and only if
\begin{equation}\label{eq:jbfqjy222}
(H+\mu^* I)x_*+c=0,~ f'_0(y_*)-\tfrac{a\mu^*}{2}=0,
\end{equation}
where $\mu^*$ is a root of the scalar function $\phi(\mu)$ defined in \eqref{eq:scalarF} in $(\max\{-\lambda_2,0\},-\lambda_1)$ such that $\phi'(\mu^*)>0$.
\end{theorem}

\emph{Proof.}
(i) We first observe that the scalar function $\phi(\mu)$ in \eqref{eq:scalarF} has the same roots as
$$
p(\mu)=\ln\varphi(\mu)-\ln\psi(\mu), ~\mu\in (\max\{-\lambda_2,0\},-\lambda_1).
$$
By Theorem \ref{th:local_pneq0}, it is necessary to assume $g_1\neq 0$. Using chain rule, \eqref{eq:phi1ds} gives
$$
p''(\mu)=\frac{\left[\sum_{i=1}^n\frac{6g_i^2}{(\lambda_i+\mu)^4}\right]\left[\sum_{i=1}^n\frac{g_i^2}{(\lambda_i+\mu)^2}\right]-{\left[\sum_{i=1}^n\frac{2g_i^2}{{(\lambda_i+\mu)}^3}\right]}^2}{[\varphi(\mu)]^2}-{\left(\ln\psi(\mu)\right)}''. $$
Define two vectors in $\R^n$:
$$
s=\bmatrix\frac{\sqrt{6}g_1}{{(\lambda_1+\mu)}^2}, \cdots, \frac{\sqrt{6}g_n}{{(\lambda_n+\mu)}^2}\endbmatrix^T, ~ t=\bmatrix\frac{g_1}{\lambda_1+\mu}, \cdots, \frac{g_n}{\lambda_n+\mu}\endbmatrix^T$$
It follows from Cauchy-Schwartz inequality that
%\begin{equation*}
%\begin{aligned}
%{\left[\sum_{i=1}^n\frac{2g_i^2}{{(\lambda_i+\mu)}^3}\right]}^2&<{(s^Tt)}^2\\
%&\le(s^Ts)(t^Tt)\\
%&=\left[\sum_{i=1}^n\frac{6g_i^2}{{(\lambda_i+\mu)}^4}\right]\left[\sum_{i=1}^n\frac{g_i^2}{{(\lambda_i+\mu)}^2}\right].
%\end{aligned}
%\end{equation*}
\[
{\left[\sum_{i=1}^n\frac{2g_i^2}{{(\lambda_i+\mu)}^3}\right]}^2
<{(s^Tt)}^2\le(s^Ts)(t^Tt)=\left[\sum_{i=1}^n\frac{6g_i^2}{{(\lambda_i+\mu)}^4}\right]\left[\sum_{i=1}^n\frac{g_i^2}{{(\lambda_i+\mu)}^2}\right].
\]
The assumption that $\psi(\mu)$ is log-concave implies that ${\left(\ln\psi(\mu)\right)}''\le 0$. Therefore, we have $p''(\mu)>0$  and hence $p(\mu)$ is strictly convex for all $\mu\in (\max\{-\lambda_2,0\},-\lambda_1)$. Thus the equation $p(\mu)=0$, as well as $\phi(\mu)=0$, has at most two real roots in the above interval, denoted by $\mu_1<\mu_2$. Suppose $\phi'(\mu_1)\ge 0$ and  $\phi'(\mu_2)\ge0$. Then, for any sufficiently small $\epsilon\in(0,(\mu_2-\mu_1)/2)$, we have
$$\phi(\mu_1+\epsilon)\ge \phi(\mu_1)=0,~ \phi(\mu_2-\epsilon)\le \phi(\mu_2)=0.$$
Therefore, there is a $\tilde{\mu}\in[\mu_1+\epsilon,\mu_2-\epsilon]$ such that $\phi(\tilde{\mu})= 0$, which is a contradiction.
Consequently, the scalar function $\phi(\mu)$ has at most one real root satisfying $\phi'(\mu)\ge 0$. Following Theorem \ref{th:local_nephi}, the proof is complete.

%\begin{remark}
%Theorem \ref{local_onlyone} is an extension of  Theorem 3.3 presented in \cite{Xia17}, where $f_0(y)=\tfrac{\sigma}{p}y^{\frac{p}{2}}$ and $\psi(\mu)=\left
%(\frac{\mu}{\sigma}\right)^{\frac{2}{p-2}}$ can be verified to be log-concave.
%\end{remark}

(ii)  Firstly, we prove that the assumption that $\psi(\mu)$ is log-concave in $(\max\{-\lambda_2,0\},-\lambda_1)$ is equivalent to
\begin{align}
f_0'''(y(\mu))+\frac{a}{ay(\mu)+b}f_0''(y(\mu))\ge0, \mu\in (\max\{-\lambda_2,0\},-\lambda_1).\label{eq:as7}
\end{align}
Applying the inverse function theorem to \eqref{eq:zjbfqj-y}, one can verify that
$$
\begin{array}{ll}
y'(\mu)=\frac{a}{2f''_0(y(\mu))},~ y''(\mu)=-\frac{af'''_0(y(\mu))y'(\mu)}{2(f''_0(y(\mu)))^2}.
\end{array}
$$
Combining with \eqref{eq:zjbfqj-psi},
$$
\begin{array}{ll}
%(\ln(\psi(\mu)))'&=\frac{\psi'(\mu)}{\psi(\mu)}=\frac{a}{2y(\mu)f''_0(y(\mu))},\\
(\ln(\psi(\mu)))''&=\frac{\psi''(\mu)\psi(\mu)-{(\psi'(\mu))}^2}{{(\psi(\mu))}^2}\\
&=-\frac{a^3}{4(ay(\mu)+b){[f_0''(y(\mu))]}^3}[f_0'''(y(\mu))+\frac{a}{ay(\mu)+b}f_0''(y(\mu))].
\end{array}
$$
which implies \eqref{eq:as7}.

It is sufficient to prove that the optimality condition in Theorem \ref{th:local_suphi} is necessary, if \eqref{eq:as7} hold. Let $(x_*,y_*)$ be a local non-global minimizer of \eqref{eq:TRS-L}. It follows from Theorem \ref{th:local_nephi} and the discussion after Theorem \ref{th:local_suphi} that \eqref{eq:jbfqjy222} holds and $\mu^*$ is a root of the scalar function $\phi$ in $(\max\{-\lambda_2,0\},-\lambda_1)$ such that $\phi'(\mu^*)\ge 0$. The remaining part is to show  $\phi'(\mu^*)>0$. Suppose this is not true, that is, we assume $\phi'(\mu^*)=0$. According to \eqref{eq:jbfqj4}, the reduced Hessian
\begin{equation*}
B:=W^TG W
\end{equation*}
has a zero eigenvalue, where
\begin{equation}\label{eq:GW}
G=\bmatrix H+\mu^*I & 0\\ 0 & f_0''(y_*)\endbmatrix,~
W=\bmatrix V & 0\\ 0 & 1\endbmatrix
\bmatrix
    \frac{g_2}{\lambda_2+\mu^*} & \frac{g_3}{\lambda_3+\mu^*} & \cdots & \frac{g_n}{\lambda_n+\mu^*}
     & \frac{a}{2}\\
     -\frac{g_1}{\lambda_1+\mu^*} & 0 & \cdots & 0
     & 0 \\
     0 & -\frac{g_1}{\lambda_1+\mu^*} & \cdots & 0
     & 0 \\
     \vdots & \vdots &   & \vdots
     & \vdots \\
     0 & 0 & \cdots & -\frac{g_1}{\lambda_1+\mu^*}
     & 0 \\
     0 & 0 & \cdots & 0  & -\frac{g_1}{\lambda_1+\mu^*}
\endbmatrix.
\end{equation}
Let $q\neq 0$ be an eigenvector of $B$ corresponding to the zero eigenvalue. That is,
\begin{equation}\label{eq:B(mu)q=0}
W^TGWq=0.
\end{equation}
Since columns of $W$ form a basis of the hyperplane $x_*^Ts-\tfrac{a}{2}t=0$, i.e.,
\begin{equation*}
W^T\bmatrix x_*^T & -\tfrac{a}{2}\endbmatrix^T=0,
\end{equation*}
and $W$ is of full column rank. It follows from \eqref{eq:B(mu)q=0} that
\begin{equation}\label{eq:GWq}
GWq=\gamma
\bmatrix
  x_*\\
   -\tfrac{a}{2}
  \endbmatrix
\end{equation}
for some $\gamma\in\R$. It follows from the linearly independence of columns of $W$ and $q\neq 0$ that
\begin{equation}\label{eq:Wqneq0}
Wq\neq 0.
\end{equation}
Notice that the matrix $G$ in \eqref{eq:GW} is nonsingular. By combining \eqref{eq:GWq} and \eqref{eq:Wqneq0}, we have
\begin{equation}\label{eq:Wqst}
Wq =\gamma G^{-1}\bmatrix
  x_*\\
   -\tfrac{a}{2}
  \endbmatrix=:\bmatrix s_*\\  t_*\endbmatrix,
\end{equation}
where $s_*\in \R^n, t_*\in \R$. Thus $\gamma\neq 0$. Furthermore, we claim that $s_*\neq 0, t_*\neq 0$ and
\begin{equation}\label{eq:st}
t_*=\frac{2}{a}x_*^Ts_*.
\end{equation}
Actually, substituting the matrix $G$ defined in \eqref{eq:GW} into the second equality of \eqref{eq:Wqst} yields  $t_*=-\tfrac{\gamma a}{2f''_0(y)}\neq 0$.
Multiplying $q^TW^T$ from left to both sides of \eqref{eq:GWq} gives
\begin{equation*}\label{eq:stgx}
0=q^TW^TGWq=\gamma(s_*^Tx_*-\tfrac{at_*}{2}),
\end{equation*}
where the first equality holds from \eqref{eq:B(mu)q=0}. Now we obtain \eqref{eq:st}. Since $t_*\neq 0$, it follows from \eqref{eq:st} that $s_*\neq 0$.
Define
$$
y(\beta):=\tfrac{\|x_*+ \beta s_*\|^2-b}{a},~ h(\beta):=q(x_*+\beta s_*)+f_0(y(\beta)).
$$
One can verify that
$$
\begin{array}{ll}
h'(\beta)&= s_*^T\nabla q(x_*+\beta s_*)+f'_0(y(\beta))y'(\beta),\\
h''(\beta)&= s_*^T\nabla^2 q(x_*+\beta s_*)s_*+f''_0(y(\beta))[y'(\beta)]^2+f'_0(y(\beta))y''(\beta),\\
h'''(\beta)&= f'''_0(y(\beta))[y'(\beta)]^3+3f''_0(y(\beta))y'(\beta)y''(\beta)+f'_0(y(\beta))y'''(\beta),
\end{array}
$$
and $y(0)=y_*, y'(0)=t_*, y''(0)=\tfrac{2s_*^Ts_*}{a}, y'''(0)=0$.
Let $\beta=0$,
$$
\begin{array}{ll}
h'(0)&= s_*^T(Hx_*+c)+f'_0(y_*)t_*,\\
h''(0)&= s_*^THs_*+f''_0(y_*)t_*^2+f'_0(y_*)\tfrac{2s_*^Ts_*}{a}=q^TW^TGWq,\\
h'''(0)&= f'''_0(y_*)t_*^3+3f''_0(y_*)t_*\tfrac{2s_*^Ts_*}{a},
\end{array}
$$
where the second equality on $h''(0)$ follows from \eqref{eq:jbfqjy222},  \eqref{eq:GW} and \eqref{eq:Wqst}.
The first-order necessary optimality condition \eqref{eq:jbfqjy222} and \eqref{eq:st} implies $h'(0)=0$. According to the definition of $q$ in \eqref{eq:B(mu)q=0}, we have $h''(0)=0$. Notice that $t_*\neq 0$. Substituting \eqref{eq:st} and $x_*^Tx_*=ay_*+b$ into $h'''(0)$ yields that
\[
\frac{h'''(0)}{t_*^3}=
f'''_0(y_*)+6f''_0(y_*)\frac{s_*^Ts_*}{at_*^2}
=f'''_0(y_*)+\frac{3af''_0(y_*)}{2(ay_*+b)}\frac{s_*^Ts_*x_*^Tx_*}{(x_*^Ts_*)^2}
\ge
f'''_0(y_*)+\frac{3af''_0(y_*)}{2(ay_*+b)}>0,
\]
where the last two inequalities follow from Cauchy-Schwartz inequality,
\eqref{eq:as7} and $\tfrac{a}{ay(\mu)+b}f''_0(y(\mu))>0$. Note that the first inequality also needs the fact $a>0$. We obtain $h'''(0)\neq 0$, which contradicts the facts that $(x_*+\beta s_*, y(\beta))$ is feasible for \eqref{eq:TRS-L} for all $\beta\in\R$, and $(x_*,y_*)$ is a local minimizer of \eqref{eq:TRS-L}.  Therefore, $\phi'(\mu^*)>0$, the proof of the necessary part is complete.
\eproof

\begin{figure}[h]
	\centering
	\subfloat[{Example \ref{ex:l1}}]{\label{logppinexample1}
		\includegraphics[width=0.45\textwidth]{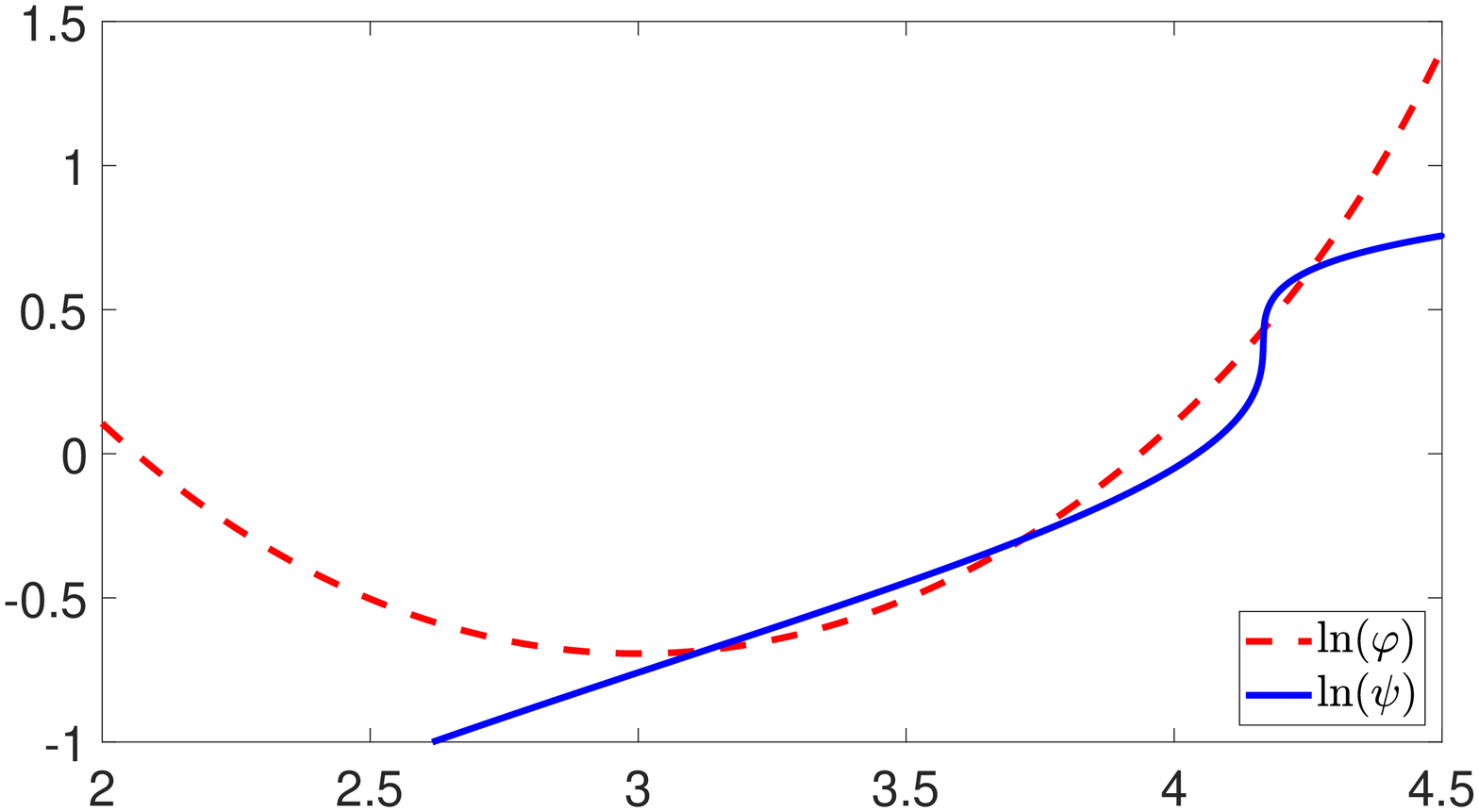}}
	~
	\subfloat[{Example \ref{ex:l2}}]{\label{logppinexample2}
		\includegraphics[width=0.45\textwidth]{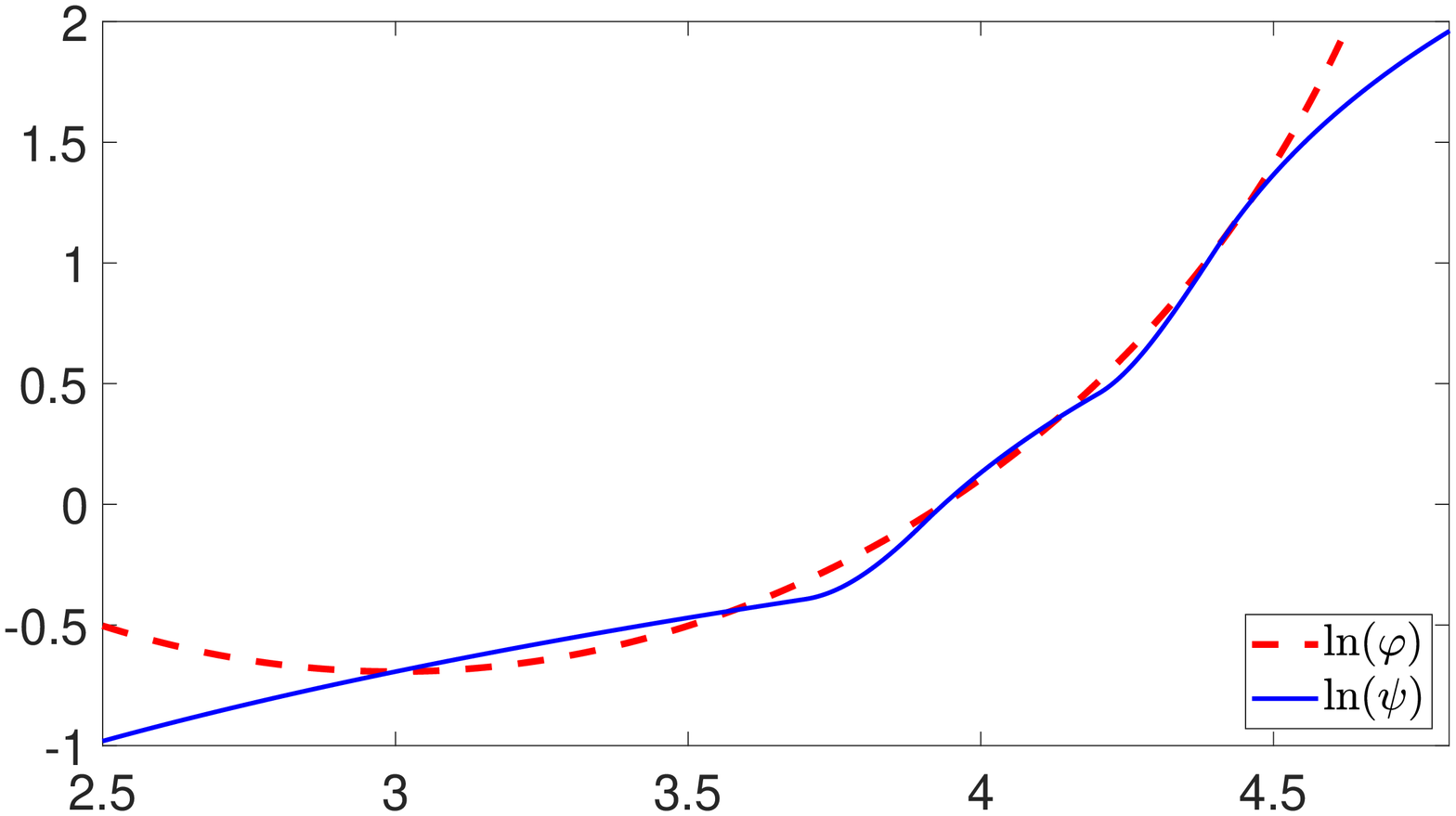}}
	\caption{Variants of $\ln\varphi(\mu)$ and $\ln\psi(\mu)$ in Examples \ref{ex:l1} and \ref{ex:l2}.}
\end{figure}

Figures \eqref{logppinexample1} and \eqref{logppinexample2} illustrate the functions $\ln\psi(\mu)$ and $\ln\varphi(\mu)$ for Examples \ref{ex:l1} and  \ref{ex:l2}, respectively. One can observe that in both cases $\ln\psi(\mu)$ is not concave, which implies the necessity of the log-concavity assumption in Theorem \ref{local_onlyone}.

Finally, we present some examples of $f_0(y)$ satisfying the assumptions in Theorem \ref{local_onlyone}.
\begin{corollary}\label{cor2}
Problem \eqref{eq:TRS-L} has at most one local non-global minimizer with necessary and sufficient optimality condition for the local non-global minimizer, if one of following conditions holds:
	\begin{enumerate}[(i)]
		\item $f_0(y)=\alpha y^d, ~\alpha>0, ~d>1, ~ a=1, ~ b=0$, (which includes Problem \eqref{eq:pRS-gtrs} as a special case).
		\item $f_0(y)$ is a strongly convex quadratic function,  $a>0, ~b\ge0$.
		\item  $f_0(y)$ is a strongly convex cubic polynomial function in $(0,+\infty)$, $a=1, ~b=0$.
	\end{enumerate}
\end{corollary}
\emph{Proof.}
According to Theorem \ref{local_onlyone}, we only need to verify that \eqref{eq:as7} holds for each case.

(i)	$f_0(y)=\alpha y^d,
\ f'_0(y)=\alpha dy^{d-1}, \ f''_0(y)=\alpha d(d-1)y^{d-2},\   f'''_0(y)=\alpha d(d-1)(d-2)y^{d-3}.$
%Let $f'_0(y)=\frac{\mu}{2}$, $y(\mu)={\left(\frac{\mu }{2\alpha d}\right)}^{\frac{1}{d-1}}$. Then
%$$\ln\psi(\mu)=\ln y(\mu)=\frac{1}{d-1}(\ln\mu-\ln(2\alpha d)),$$
%which is concave on $(0, +\infty)$ for $d>1$. The assumption of Theorem \ref{local_onlyone} is satisfied. Moreover,
	$$ f'''_0(y)+\frac{1}{y}f''_0(y)=\alpha d{(d-1)}^2y^{d-3}>0,\ \forall y\in(0,+\infty),$$
	which implies that the assumption of Theorem  \ref{local_onlyone} holds.

(ii) Without loss of generality, we assume $f_0(y)=\alpha y^2+\beta y, ~\alpha>0$. Then  $f'_0(y)=2\alpha y+\beta, f''_0(y)=2\alpha, f'''_0(y)=0$.
%	Let $f'_0(y)=\frac{a\mu}{2}$, $y(\mu)=\frac{a\mu/2-\beta}{2\alpha}$. Thus,  $$\psi(\mu)=ay(\mu)+b=\frac{a^2}{4\alpha}\mu-\frac{a\beta}{2\alpha}+b,$$
%	which is an increasing linear function in terms of $\mu$  and positive as $y(\mu)\in S$. Thus, $\psi(\mu)$ is log-concave.
%Moreover, we have
$$ f'''_0(y)+\frac{a}{ay+b}f''_0(y)=\frac{2a\alpha}{ay+b}>0,~ \forall y\in (-\tfrac{b}{a}, +\infty).$$
	
(iii) Assume that $f_0(y)=\alpha y^3+\beta y^2+\gamma y$ is strongly convex on $(0,+\infty)$. Notice that
$$
f'_0(y)=3\alpha y^2+2\beta y+\gamma,  \ f''_0(y)=6\alpha y+2\beta, \ f'''_0(y)=6\alpha.
$$
Since $f''_0(y)=6\alpha \left(y+\frac{\beta}{3\alpha }\right)>0, ~\forall y\in (0,+\infty)$. It turns out that $\alpha>0$ and  $\beta>0$.
	
%Let $f'_0(y)=\frac{\mu}{2}$. Then we have \[
%(\mu)=\frac{\sqrt{\beta^2+3\alpha(\frac{\mu}{2}-\gamma)}}{3\alpha}-\frac{\beta}{3\alpha%},
%\]
%which implies that
%$$\ln\psi(\mu)=\ln y(\mu)=\ln\left(\sqrt{\mu+\frac{2\beta^2}{3\alpha}-2\gamma}-\sqrt{\frac{2\beta^2}{3\alpha}}
%	\right)-\frac{\ln(3\alpha)+\ln{2}}{2}.$$
%It can be verified that $\ln\psi(\mu)$ is concave in $R$. Moreover, we have

$$f'''_0(y)+\frac{1}{y}f''_0(y)>0+0=0, ~\forall y\in (0,+\infty).$$
\eproof

Note that we cannot extend the above cases presented in Corollary \ref{cor2} to the general quartic polynomial case, see
Example \ref{ex:l1} for a counterexample.

\section{Conclusion}
We raise a fundamental question (Question \ref{ques}) whether  local optimality can be checked in polynomial time for hidden convex optimization.  Then we focus on the newly proposed optimization
problem by jointing nonconvex trust-region subproblem with  convex optimization \eqref{eq:TRS-C}.
We present a necessary and sufficient optimality condition for global minimizer of \eqref{eq:TRS-C}, which reveals its hidden convexity. However, it is difficult to establish a necessary and sufficient optimality condition for
local non-global minimizer, except for some quadratic single-constraint cases.
Moreover, different from  trust region subproblem (which has at most one local non-global minimizer) and convex optimization (without local non-global minimizer), their joint problem could have more than one local non-global minimizer. There is a quartic polynomial case of \eqref{eq:TRS-C} with two local non-global minimizers.  We then present a general approach to generate the instances with arbitrary number of local non-global minimizers. Consequently, we conclude that the existence of many local minimizers
is NOT a (or at least not a unique) reason making global optimization difficult.

While we have present some negative  evidences, Question \ref{ques} remains open, even on the joint problem problem \eqref{eq:TRS-C}.

%\begin{enumerate}[\textbf{open question} 1]
%\item{Does problem \eqref{eq:TRS-L} has sufficient and necessary condition for its local non-global minimizer?}
%%\item{Can the conditions of Theorem \ref{local_onlyone} or \ref{local_sune} be further relaxed?}
%
%\item{Complexity of testing whether a given point is a local minimizer of \eqref{eq:TRS-C}.}
%
%%\item{Can the assumptions that $f_i$ is twice continuously differential at $y_*$ for $i=0,1,\cdots,m$ in Theorem \ref{th:local_peq0} be relaxed or removed?}
%\end{enumerate}
% Appendix here
% Options are (1) APPENDIX (with or without general title) or
%             (2) APPENDICES (if it has more than one unrelated sections)
% Outcomment the appropriate case if necessary
%
% \begin{APPENDIX}{<Title of the Appendix>}
% \end{APPENDIX}
%
%   or
%
% \begin{APPENDICES}
% \section{<Title of Section A>}
% \section{<Title of Section B>}
% etc
% \end{APPENDICES}

% Acknowledgments here
\section*{Acknowledgments.}\label{se:7}
% Enter the text of acknowledgments here
This research was supported by the National Natural Science Foundation of China under grants 11822103, 11571029, and 11771056, and by the Beijing Natural Science Foundation, grant Z180005.
% References here (outcomment the appropriate case)

% CASE 1: BiBTeX used to constantly update the references
%   (while the paper is being written).
\bibliographystyle{informs2014} % outcomment this and next line in Case 1
\bibliography{reference} % if more than one, comma separated

% CASE 2: BiBTeX used to generate mypaper.bbl (to be further fine tuned)
%\input{mypaper.bbl} % outcomment this line in Case 2

\end{document}